\documentclass[10pt]{article}
\usepackage{setspace}

\usepackage{amsmath} 
\usepackage{amssymb}
\usepackage{latexsym}
\usepackage{xcolor}
\usepackage{fancyhdr}
\allowdisplaybreaks 

 \usepackage{comment}

\def\XXint#1#2#3{{\setbox0=\hbox{$#1{#2#3}{\int}$} 
\vcenter{\hbox{$#2#3$}}\kern-.5\wd0}}   

 \numberwithin{equation}{section}
\newtheorem{theorem}[equation]{Theorem}
\newtheorem{proposition}[equation]{Proposition}
\newtheorem{definition}[equation]{Definition}
\newtheorem{remark}[equation]{Remark}
\newtheorem{lemma}[equation]{Lemma}

\title{
A  uniqueness theorem  for nonvariational solutions of the
Helmholtz equation} 
 
\author{  
Massimo Lanza de Cristoforis
\\
Dipartimento di Matematica `Tullio Levi-Civita', 
\\
Universit\`a degli Studi di Padova, 
\\
Via Trieste 63, Padova 35121, 
Italy. 
\\
E-mail: mldc@math.unipd.it   }

\date{\ }

\begin{document}
 
 \maketitle

\noindent
{\bf Abstract:}  We consider   a bounded open subset $\Omega$ of ${\mathbb{R}}^n$ of class $C^{1,\alpha}$ for some 
$\alpha\in]0,1[$, and we define a distributional outward unit normal derivative for 
 $\alpha$-H\"{o}lder continuous solutions of the Helmholtz equation in the exterior  of $\Omega$ that may   not have a classical outward unit normal derivative at the boundary points of $\Omega$ and that may have 
an infinite Dirichlet integral around the boundary of $\Omega$. Namely for solutions  that  do not belong to the classical variational setting.

Then we show
a Schauder boundary regularity result for $\alpha$-H\"{o}lder continuous functions that have the Laplace operator in a Schauder space of negative exponent and we prove  a uniqueness theorem for $\alpha$-H\"{o}lder continuous solutions 
of
the exterior Dirichlet and impedance boundary value problems for the Helmholtz equation
that satisfy the Sommerfeld radiation condition at infinity
 in the above mentioned nonvariational setting.   
 \vspace{\baselineskip}

 

\section{Introduction} Our starting point are the classical examples of Prym \cite{Pr1871} and  Hadamard \cite{Ha1906} of harmonic functions in the unit ball of the plane that are continuous up to the boundary and have  infinite Dirichlet integral, \textit{i.e.},  whose
 gradient is not square-summable. Such functions solve the classical Dirichlet problem in the unit ball, but not the corresponding weak (variational) formulation. For a discussion on this point we refer to 
Maz’ya and Shaposnikova \cite{MaSh98}, Bottazzini and Gray \cite{BoGr13} and
 Bramati,  Dalla Riva and  Luczak~\cite{BrDaLu23} which contains examples of H\"{o}lder continuous harmonic functions with infinite Dirichlet integral.

In the papers \cite{La24b},\cite{La24c}, the author has analyzed the Neumann problem for the Laplace and the Poisson equation for $\alpha$-H\"{o}lder continuous solutions that may have infinite Dirichlet integral and
the present paper is a first step in the corresponding analysis  of   the Helmholtz equation.

We consider a bounded open subset $\Omega$ of ${\mathbb{R}}^n$ of class $C^{1,\alpha}$,  its exterior
\[
\Omega^-\equiv {\mathbb{R}}^n\setminus\overline{\Omega}\,,
\] 
 the space $C^{-1,\alpha}_{{\mathrm{loc}}	}(\overline{\Omega^-})$
  of sums of locally $\alpha$-H\"{o}lder continuous functions in $\overline{\Omega^-}$ and of first order partial distributional derivatives in $\Omega^-$ of   locally $\alpha$-H\"{o}lder continuous functions in $\overline{\Omega^-}$. Then we introduce a   distributional normal derivative $\partial_{\nu_{\Omega^-}}$ on $\partial\Omega$ for functions $u$ in the space $C^{0,\alpha}_{{\mathrm{loc}}	}(\overline{\Omega^-})_\Delta$ of functions in the space $C^{0,\alpha}_{{\mathrm{loc}}	}(\overline{\Omega^-})$ of locally $\alpha$-H\"{o}lder continuous functions in $\overline{\Omega^-}$ such that the distributional Laplace operator 
  $\Delta u$ in $\Omega^-$ belongs to $C^{-1,\alpha}_{{\mathrm{loc}}	}(\overline{\Omega^-})$ by extending an argument of 
    \cite[\S 5]{La24c} (see Definition \ref{defn:endedr}) and we prove the existence of a right inverse for $\partial_{\nu_{\Omega^-}}$ (see Proposition \ref{prop:recodnu}).
 
  Then we prove that a function  of $C^{0,\alpha}_{{\mathrm{loc}}	}(\overline{\Omega^-})_\Delta$ whose outward distributional normal derivative $\partial_{\nu_{\Omega^-}}$ 
  is of class $C^{0,\alpha}(\partial\Omega)$ is actually of class 
$C^{1,\alpha}_{{\mathrm{loc}}	}(\overline{\Omega^-})$ and that a corresponding result holds in $\Omega$ (cf. the Schauder regularity Theorem \ref{thm:schreim}). 

Finally, we extend to solutions of class $C^{0,\alpha}_{{\mathrm{loc}}	}(\overline{\Omega^-})_\Delta$ of the Helmholtz equation with impedance and Dirichlet boundary conditions a classical uniqueness theorem that is known to hold in case the outward unit normal derivative of the solutions  is known either  to  exist classically on the points of $\partial\Omega$ or to exist in a weak  form  in case the Dirichlet integral  of the solutions  close to    $\partial\Omega$ is finite, as it happens in the variational setting (cf.~\textit{e.g.}, Colton and Kress~\cite[Thms.~3.13, 3.37]{CoKr13}, Neittaanm\"{a}ki,   Roach~\cite{NiRo87}, N\'{e}d\'{e}lec~\cite[Thm.~2.6.5, p.~108]{Ne01}).

Instead, we plan  to develop the consequent  treatment of the existence of solutions for the above mentioned boundary value problems for the Helmholtz equations   in forthcoming papers (cf.~\cite{La25a}, \cite{La25b}). 
 
The paper is organized as follows. Section \ref{sec:prelnot} is a section of
 preliminaries and notation. 
  In section \ref{sec:doun} we introduce a distributional form of the outward unit normal derivative for functions of $C^{0,\alpha}_{{\mathrm{loc}}	}(\overline{\Omega^-})_\Delta$. In section \ref{sec:nvschre}, we prove   the Schauder regularity Theorem \ref{thm:schreim}. Section \ref{sec:racouh} is devoted to the uniqueness theorem of the Helmholtz equation with  Dirichlet and impedance boundary conditions. In the appendix at the end of the paper, we have included a formulation of the classical uniqueness theorem for the Helmholtz equation that we need in the paper.

\section{Preliminaries and notation}\label{sec:prelnot} Unless otherwise specified,  we assume  throughout the paper that
\[
n\in {\mathbb{N}}\setminus\{0,1\}\,,
\]
where ${\mathbb{N}}$ denotes the set of natural numbers including $0$. 
If $X$ and $Y$, $Z$ are normed spaces, then ${\mathcal{L}}(X,Y)$ denotes the space of linear and continuous maps from $X$ to $Y$ and ${\mathcal{L}}^{(2)}(X\times Y, Z)$ 
denotes the space of bilinear and continuous maps from $X\times Y$ to $Z$ with their usual operator norm (cf.~\textit{e.g.}, \cite[pp.~16, 621]{DaLaMu21}).

$|A|$ denotes the operator norm of a matrix $A$ with real (or complex) entries, 
       $A^{t}$ denotes the transpose matrix of $A$.	 
 Let $\Omega$ be an open subset of ${\mathbb{R}}^n$. $C^{1}(\Omega)$ denotes the set of continuously differentiable functions from $\Omega$ to ${\mathbb{C}}$. 
 Let $s\in {\mathbb{N}}\setminus\{0\}$, $f\in \left(C^{1}(\Omega)\right)^{s} $. Then   $Df$ denotes the Jacobian matrix of $f$.   
 
 For the (classical) definition of   open Lipschitz subset of ${\mathbb{R}}^n$ and of   open subset of ${\mathbb{R}}^n$ 
   of class $C^{m}$ or of class $C^{m,\alpha}$
  and of the H\"{o}lder and Schauder spaces $C^{m,\alpha}(\overline{\Omega})$
  on the closure $\overline{\Omega}$ of  an open set $\Omega$ and 
  of the H\"{o}lder and Schauder spaces
   $C^{m,\alpha}(\partial\Omega)$ 
on the boundary $\partial\Omega$ of an open set $\Omega$ for some $m\in{\mathbb{N}}$, $\alpha\in ]0,1]$, we refer for example to
    Dalla Riva, the author and Musolino  \cite[\S 2.3, \S 2.6, \S 2.7, \S 2.9, \S 2.11, \S 2.13,   \S 2.20]{DaLaMu21}.  If $m\in {\mathbb{N}}$, 
 $C^{m}_b(\overline{\Omega})$ denotes the space of $m$-times continuously differentiable functions from $\Omega$ to ${\mathbb{C}}$ such that all 
the partial derivatives up to order $m$ have a bounded continuous extension to    $\overline{\Omega}$ and we set
\[
\|f\|_{   C^{m}_{b}(
\overline{\Omega} )   }\equiv
\sum_{|\eta|\leq m}\, \sup_{x\in \overline{\Omega}}|D^{\eta}f(x)|
\qquad\forall f\in C^{m}_{b}(
\overline{\Omega} )\,.
\]
If $\alpha\in ]0,1]$, then 
$C^{m,\alpha}_b(\overline{\Omega})$ denotes the space of functions of $C^{m}_{b}(
\overline{\Omega}) $  such that the  partial derivatives of order $m$ are $\alpha$-H\"{o}lder continuous in $\Omega$. Then we equip $C^{m,\alpha}_{b}(\overline{\Omega})$ with the norm
\[
\|f\|_{  C^{m,\alpha}_{b}(\overline{\Omega})  }\equiv 
\|f\|_{  C^{m }_{b}(\overline{\Omega})  }
+\sum_{|\eta|=m}|D^{\eta}f|_{\alpha}\qquad\forall f\in C^{m,\alpha}_{b}(\overline{\Omega})\,,
\]
where $|D^{\eta}f|_{\alpha}$ denotes the $\alpha$-H\"{o}lder constant of the partial derivative $D^{\eta}f$ of order $\eta$ of $f$ in $\Omega$. If $\Omega$ is bounded, we obviously have $C^{m }_{b}(\overline{\Omega})=C^{m } (\overline{\Omega})$ and $C^{m,\alpha}_{b}(\overline{\Omega})=C^{m,\alpha} (\overline{\Omega})$. 
  Then $C^{m,\alpha}_{{\mathrm{loc}}}(\overline{\Omega }) $\index{$C^{m,\alpha}_{{\mathrm{loc}}}(\overline{\Omega }) $} denotes 
the space  of those functions $f\in C^{m}(\overline{\Omega} ) $ such that $f_{|\overline{ \Omega\cap{\mathbb{B}}_{n}(0,\rho) }} $ belongs to $
C^{m,\alpha}(   \overline{ \Omega\cap{\mathbb{B}}_{n}(0,\rho) }   )$ for all $\rho\in]0,+\infty[$.
The space of complex valued functions of class $C^m$ with compact support in an open set $\Omega$ of ${\mathbb{R}}^n$ is denoted $C^m_c(\Omega)$ and similarly for $C^\infty_c(\Omega)$. We also set ${\mathcal{D}}(\Omega)\equiv C^\infty_c(\Omega)$. Then the dual ${\mathcal{D}}'(\Omega)$ is known to be the space of distributions in $\Omega$. The support of a function or of a distribution is denoted by the abbreviation `${\mathrm{supp}}$'.  We also set
\[
\Omega^-\equiv {\mathbb{R}}^n\setminus\overline{\Omega}\,.
\] 
We denote by $\nu_\Omega$ or simply by $\nu$ the outward unit normal of $\Omega$ on $\partial\Omega$. Then $\nu_{\Omega^-}=-\nu_\Omega$ is the outward unit normal of $\Omega^-$ on $\partial\Omega=\partial\Omega^-$.
We now summarize the definition and some elementary properties of the Schau\-der space $C^{-1,\alpha}(\overline{\Omega})$ 
by following the presentation of Dalla Riva, the author and Musolino \cite[\S 2.22]{DaLaMu21}.
\begin{definition} 
\label{defn:sch-1}\index{Schauder space!with negative exponent}
 Let $\alpha\in]0,1]$. Let $\Omega$ be a bounded open subset of ${\mathbb{R}}^{n}$. We denote by $C^{-1,\alpha}(\overline{\Omega})$ the subspace 
 \[
 \left\{
 f_{0}+\sum_{j=1}^{n}\frac{\partial}{\partial x_{j}}f_{j}:\,f_{j}\in 
 C^{0,\alpha}(\overline{\Omega})\ \forall j\in\{0,\dots,n\}
 \right\}\,,
 \]
 of the space of distributions ${\mathcal{D}}'(\Omega)$  in $\Omega$ and we set
 \begin{eqnarray}
\label{defn:sch-2}
\lefteqn{
\|f\|_{  C^{-1,\alpha}(\overline{\Omega})  }
\equiv\inf\biggl\{\biggr.
\sum_{j=0}^{n}\|f_{j}\|_{ C^{0,\alpha}(\overline{\Omega})  }
:\,
}
\\ \nonumber
&&\qquad\qquad\qquad\qquad
f=f_{0}+\sum_{j=1}^{n}\frac{\partial}{\partial x_{j}}f_{j}\,,\ 
f_{j}\in C^{0,\alpha}(\overline{\Omega})\ \forall j\in \{0,\dots,n\}
\biggl.\biggr\}\,.
\end{eqnarray}
\end{definition}
$(C^{-1,\alpha}(\overline{\Omega}), \|\cdot\|_{  C^{-1,\alpha}(\overline{\Omega})  })$ is known to be a Banach space and  is continuously embedded into ${\mathcal{D}}'(\Omega)$. Also, the definition of the norm $\|\cdot\|_{  C^{-1,\alpha}(\overline{\Omega})  }$ implies that $C^{0,\alpha}(\overline{\Omega})$ is continuously embedded into $C^{-1,\alpha}(\overline{\Omega})$ and that the partial differentiation $\frac{\partial}{\partial x_{j}}$ is continuous from 
$C^{0,\alpha}(\overline{\Omega})$ to $C^{-1,\alpha}(\overline{\Omega})$ for all $j\in\{1,\dots,n\}$.  Generically, the  elements of $C^{-1,\alpha}(\overline{\Omega})$ are not integrable functions, but distributions in $\Omega$.  Then we have the following statement of \cite[Prop.~3.1]{La24c} that shows that the elements of $C^{-1,\alpha}(\overline{\Omega}) $, which belong to the dual of ${\mathcal{D}}(\Omega)$, can be extended to elements of the dual of $C^{1,\alpha}(\overline{\Omega})$.  
\begin{proposition}\label{prop:nschext}
 Let $\alpha\in]0,1]$. Let $\Omega$ be a bounded open Lipschitz subset of ${\mathbb{R}}^{n}$.  There exists one and only one  linear and continuous extension operator $E^\sharp_\Omega$ from $C^{-1,\alpha}(\overline{\Omega})$ to $\left(C^{1,\alpha}(\overline{\Omega})\right)'$ such that
 \begin{eqnarray}\label{prop:nschext2}
\lefteqn{
\langle E^\sharp_\Omega[f],v\rangle
}
\\ \nonumber
&&\ \
 =
\int_{\Omega}f_{0}v\,dx+\int_{\partial\Omega}\sum_{j=1}^{n} (\nu_{\Omega})_{j}f_{j}v\,d\sigma
 -\sum_{j=1}^{n}\int_{\Omega}f_{j}\frac{\partial v}{\partial x_j}\,dx
\quad \forall v\in C^{1,\alpha}(\overline{\Omega})
\end{eqnarray}
for all $f=  f_{0}+\sum_{j=1}^{n}\frac{\partial}{\partial x_{j}}f_{j}\in C^{-1,\alpha}(\overline{\Omega}) $. Moreover, 
\begin{equation}\label{prop:nschext1}
E^\sharp_\Omega[f]_{|\Omega}=f\,, \ i.e.,\ 
\langle E^\sharp_\Omega[f],v\rangle=\langle f,v\rangle \qquad\forall v\in {\mathcal{D}}(\Omega)
\end{equation}
for all $f\in C^{-1,\alpha}(\overline{\Omega})$ and
\begin{equation}\label{prop:nschext3}
\langle E^\sharp_\Omega[f],v\rangle =\langle f,v\rangle \qquad\forall v\in C^{1,\alpha}(\overline{\Omega})
\end{equation}
for all $f\in C^{0,\alpha}(\overline{\Omega})$.
\end{proposition}
When no ambiguity can arise, we simply write $E^\sharp$ instead of $E^\sharp_\Omega$. To see why the extension operator $E^\sharp$ can be considered as `canonical', we refer to \cite[\S 2]{La24d}. Next we prove the following multiplication lemma.
\begin{lemma}\label{lem:multc1ac-1a}
  Let   $\alpha\in ]0,1[$. Let $\Omega$ be a bounded open Lipschitz subset of ${\mathbb{R}}^{n}$. Then the pointwise product is bilinear and continuous from 
  $C^{1,\alpha}(\overline{\Omega})\times C^{-1,\alpha}(\overline{\Omega})$ to $C^{-1,\alpha}(\overline{\Omega})$. 
\end{lemma}
{\bf Proof.} Let $(f,g)\in C^{1,\alpha}(\overline{\Omega})\times C^{-1,\alpha}(\overline{\Omega})$. By definition of 
$C^{-1,\alpha}(\overline{\Omega})$, there exist $g_j\in C^{0,\alpha}(\overline{\Omega})$ for $j\in\{0,1,\dots,n\}$ such that
\begin{equation}\label{lem:multc1ac-1a1}
g=g_0+\sum_{j=1}^n\frac{\partial g_j}{\partial x_j}
\end{equation}
and thus we have
\[
fg=fg_0+\sum_{j=1}^nf\frac{\partial g_j}{\partial x_j}
=fg_0+\sum_{j=1}^n \frac{\partial (fg_j)}{\partial x_j}-\sum_{j=1}^n g_j\frac{\partial f }{\partial x_j}
\]
and the membership of $fg_j$, $g_j\frac{\partial f }{\partial x_j}$ in $C^{0,\alpha}(\overline{\Omega})$ for all $j\in\{0,1,\dots,n\}$ implies that $fg\in C^{-1,\alpha}(\overline{\Omega})$. Now let $c_1$, $c_2$ denote the norm of the pointwise product in $C^{0,\alpha}(\overline{\Omega})$ and of the embedding of $C^{1,\alpha}(\overline{\Omega})$ into $C^{0,\alpha}(\overline{\Omega})$, respectively.
 Then we have 
\begin{eqnarray*}
\lefteqn{
\|fg\|_{C^{-1,\alpha}(\overline{\Omega})}
}
\\ \nonumber
&&\ \ 
\leq \left\|fg_0-\sum_{j=1}^n g_j\frac{\partial f }{\partial x_j}\right\|_{C^{0,\alpha}(\overline{\Omega})}
+\sum_{j=1}^n\left\|fg_j\right\|_{C^{0,\alpha}(\overline{\Omega})}
\\ \nonumber
&&\ \ 
\leq 2c_1\left(\|f\|_{C^{0,\alpha}(\overline{\Omega})}+\sum_{j=1}^n\left\|\frac{\partial f }{\partial x_j}\right\|_{C^{0,\alpha}(\overline{\Omega})}
\right)\left(\|g_0\|_{C^{0,\alpha}(\overline{\Omega})}+\sum_{j=1}^n\left\|g_j\right\|_{C^{0,\alpha}(\overline{\Omega})}
\right)\,.
\end{eqnarray*}
Then by taking the infimum on all possible $g_j$ as in (\ref{lem:multc1ac-1a1}), we obtain 
\[
\|fg\|_{C^{-1,\alpha}(\overline{\Omega})}\leq 2c_1(c_2+n)\|f\|_{
C^{1,\alpha}(\overline{\Omega})
}\|g\|_{C^{-1,\alpha}(\overline{\Omega})}
\]
and thus the product is continuous. \hfill  $\Box$ 

\vspace{\baselineskip}

Next we introduce our space for the solutions.

\begin{definition}\label{defn:c0ade}
 Let   $\alpha\in ]0,1]$. Let $\Omega$ be a bounded open subset of ${\mathbb{R}}^{n}$. Let
 \begin{eqnarray}\label{defn:c0ade1}
C^{0,\alpha}(\overline{\Omega})_\Delta
&\equiv&\biggl\{u\in C^{0,\alpha}(\overline{\Omega}):\,\Delta u\in C^{-1,\alpha}(\overline{\Omega})\biggr\}\,,
\\ \nonumber
\|u\|_{ C^{0,\alpha}(\overline{\Omega})_\Delta }
&\equiv& \|u\|_{ C^{0,\alpha}(\overline{\Omega})}
+\|\Delta u\|_{C^{-1,\alpha}(\overline{\Omega})}
\qquad\forall u\in C^{0,\alpha}(\overline{\Omega})_\Delta\,.
\end{eqnarray}
\end{definition}
Since $C^{0,\alpha}(\overline{\Omega})$ and $C^{-1,\alpha}(\overline{\Omega})$ are Banach spaces,   $\left(\|u\|_{ C^{0,\alpha}(\overline{\Omega})_\Delta }, \|\cdot \|_{ C^{0,\alpha}(\overline{\Omega})_\Delta }\right)$ is a Banach space if $\Omega$ is a bounded open subset of ${\mathbb{R}}^{n}$.  
For subsets $\Omega$ that are not necessarily bounded, we introduce the following statement.
\begin{lemma}\label{lem:c1alcof}
 Let   $\alpha\in ]0,1]	$. Let $\Omega$ be an open  subset of ${\mathbb{R}}^{n}$. Then the space
\begin{eqnarray}\label{lem:c1alcof1}
\lefteqn{
 C^{0,\alpha}_{	{\mathrm{loc}}	}(\overline{\Omega})_\Delta\equiv\biggl\{
 f\in C^{0}(\overline{\Omega}):\, f_{|\overline{\Omega'}}\in C^{0,\alpha}(\overline{\Omega'})_\Delta\ \text{for\ all\  } 
 }
\\ \nonumber
&&\qquad\qquad\qquad\qquad 
 \text{bounded\ open\ subsets}\ \Omega'\ \text{of}\ {\mathbb{R}}^n\ \text{such\ that} \  \overline{\Omega'}\subseteq\overline{\Omega}  \biggr\}\,,
\end{eqnarray}
 with the family of seminorms
\begin{eqnarray}\label{lem:c1alcof2}
\lefteqn{
{\mathcal{P}}_{C^{0,\alpha}_{	{\mathrm{loc}}	}(\overline{\Omega})_\Delta}\equiv
\biggl\{
\|\cdot\|_{C^{0,\alpha}(\overline{\Omega'})_\Delta}:\,  
}
\\ \nonumber
&& \qquad\qquad\qquad\qquad 
 \Omega'\ \text{is\ a\ bounded\ open  subset\ of\ }  {\mathbb{R}}^n \ \text{such\ that} \  \overline{\Omega'}\subseteq\overline{\Omega}
\biggr\}
\end{eqnarray}
is a Fr\'{e}chet space.
\end{lemma}
{\bf Proof.} Since $\overline{\Omega}$ is the union of an increasing sequence of  sets $\overline{\Omega'}$ with $\Omega'$ as in (\ref{lem:c1alcof2}), the family of seminorms in (\ref{lem:c1alcof2}) is equivalent to a countable family of seminorms as in (\ref{lem:c1alcof2}) and accordingly,  $C^{0,\alpha}_{	{\mathrm{loc}}	}(\overline{\Omega})_\Delta$ with the above family of seminorms is metrizable. In order to prove that $C^{0,\alpha}_{	{\mathrm{loc}}	}(\overline{\Omega})_\Delta$ is complete it suffices to show that each Cauchy sequence is convergent and to do so 
it suffices to observe that $C^{0,\alpha}(\overline{\Omega'})_\Delta$ is a Banach space for all $\Omega'$ as in (\ref{lem:c1alcof1}).\hfill  $\Box$ 

\vspace{\baselineskip}

Next we prove the following multiplication lemma.
\begin{lemma}\label{lem:multc2ac0ad}
  Let   $\alpha\in ]0,1[$. Let $\Omega$ be a bounded open Lipschitz subset of ${\mathbb{R}}^{n}$. Then the pointwise product is bilinear and continuous from 
  $C^{2,\alpha}(\overline{\Omega})\times C^{0,\alpha}(\overline{\Omega})_\Delta$ to $C^{0,\alpha}(\overline{\Omega})_\Delta$. 
\end{lemma}
{\bf Proof.} Since both $C^{2,\alpha}(\overline{\Omega})$ and $C^{0,\alpha}(\overline{\Omega})_\Delta$ are continuously embedded into $C^{0,\alpha}(\overline{\Omega})$ and the pointwise product is continuous in  $C^{0,\alpha}(\overline{\Omega})$, then there exists $c_0\in]0,+\infty[$ such that
\begin{equation}\label{lem:multc2ac0ad1}
\|fg\|_{C^{0,\alpha}(\overline{\Omega})}\leq c_0\|f\|_{C^{2,\alpha}(\overline{\Omega})}\|g\|_{C^{0,\alpha}(\overline{\Omega})_\Delta}\qquad\forall (f,g)\in C^{2,\alpha}(\overline{\Omega})\times C^{0,\alpha}(\overline{\Omega})_\Delta\,.
\end{equation}
If  $(f,g)\in C^{2,\alpha}(\overline{\Omega})\times C^{0,\alpha}(\overline{\Omega})_\Delta$, then we have
\[
\Delta (fg)=f\Delta g +2Df(Dg)^t+g\Delta f\,.
\] 
Moreover, the membership of $f$ in $ C^{2,\alpha}(\overline{\Omega}) $
implies that $f\in C^{1,\alpha}(\overline{\Omega})$, $Df\in (C^{1,\alpha}(\overline{\Omega}))^n$, $\Delta f\in C^{0,\alpha}(\overline{\Omega})$,  and the membership of $g$ in $C^{0,\alpha}(\overline{\Omega})_\Delta$ implies that $g\in  C^{0,\alpha}(\overline{\Omega})$, $Dg\in (C^{-1,\alpha}(\overline{\Omega}))^n$, 
 $\Delta g\in  C^{-1,\alpha}(\overline{\Omega})$ for all $(f,g)\in C^{2,\alpha}(\overline{\Omega})\times C^{0,\alpha}(\overline{\Omega})_\Delta$. Then the continuity of the pointwise product in  $C^{0,\alpha}(\overline{\Omega})$ and the product Lemma \ref{lem:multc1ac-1a} imply that $\Delta (fg)\in C^{0,\alpha}(\overline{\Omega})$ and that there exists $c'\in]0,+\infty[$ such that
\begin{eqnarray}\label{lem:multc2ac0ad2}
\lefteqn{
\|\Delta (fg)\|_{C^{-1,\alpha}(\overline{\Omega})}\leq c'
\biggl(\|f\|_{C^{1,\alpha}(\overline{\Omega})}\|\Delta g\|_{C^{-1,\alpha}(\overline{\Omega})}
}
\\ \nonumber
&&\qquad
+
2\sum_{j=1}^n\left\|\frac{\partial f}{\partial x_j}\right\|_{C^{1,\alpha}(\overline{\Omega})}\left\|\frac{\partial g}{\partial x_j}\right\|_{C^{-1,\alpha}(\overline{\Omega})}
+\|\Delta f\|_{C^{0,\alpha}(\overline{\Omega})}\|g\|_{C^{0,\alpha}(\overline{\Omega})}
\biggr)\,,
\end{eqnarray}
for all $(f,g)\in C^{2,\alpha}(\overline{\Omega})\times C^{0,\alpha}(\overline{\Omega})_\Delta$.
Then the continuity of the embedding of $ C^{2,\alpha}(\overline{\Omega})$ into $ C^{1,\alpha}(\overline{\Omega})$, 
the continuity of $\frac{\partial  }{\partial x_j}$ from  $ C^{2,\alpha}(\overline{\Omega})$ to $ C^{1,\alpha}(\overline{\Omega})$
and from $ C^{0,\alpha}(\overline{\Omega})$ to $ C^{-1,\alpha}(\overline{\Omega})$, of $\Delta$ from  $C^{2,\alpha}(\overline{\Omega})$ to $ C^{0,\alpha}(\overline{\Omega})$
and the definition of norm in $C^{0,\alpha}(\overline{\Omega})_\Delta$ imply that there exists 
$c''\in]0,+\infty[$ such that
\begin{eqnarray}\label{lem:multc2ac0ad3}
\lefteqn{
\biggl(\|f\|_{C^{1,\alpha}(\overline{\Omega})}\|\Delta g\|_{C^{-1,\alpha}(\overline{\Omega})}
}
\\ \nonumber
&&\qquad
+
2\sum_{j=1}^n\left\|\frac{\partial f}{\partial x_j}\right\|_{C^{1,\alpha}(\overline{\Omega})}\left\|\frac{\partial g}{\partial x_j}\right\|_{C^{-1,\alpha}(\overline{\Omega})}
+\|\Delta f\|_{C^{0,\alpha}(\overline{\Omega})}\|g\|_{C^{0,\alpha}(\overline{\Omega})}
\biggr)
\\ \nonumber
&&\qquad
\leq
c''\|f\|_{C^{2,\alpha}(\overline{\Omega})}
\biggl(
\| g\|_{C^{0,\alpha}(\overline{\Omega})_\Delta}+2n\| g\|_{C^{0,\alpha}(\overline{\Omega})}+\| g\|_{C^{0,\alpha}(\overline{\Omega})}
\biggr)\,,
\end{eqnarray}
for all $(f,g)\in C^{2,\alpha}(\overline{\Omega})\times C^{0,\alpha}(\overline{\Omega})_\Delta$.
Then by combining inequalities (\ref{lem:multc2ac0ad1})--(\ref{lem:multc2ac0ad3}), we deduce the continuity of the pointwise product as in the statement. \hfill  $\Box$ 

\vspace{\baselineskip}

We also note that the following elementary lemma holds.
\begin{lemma}\label{lem:c1alf}
 Let   $\alpha\in ]0,1[$. Let $\Omega$ be an open  subset of ${\mathbb{R}}^{n}$. Then the space
\begin{eqnarray}\label{lem:c1alf1}
\lefteqn{
 C^{1,\alpha}_{	{\mathrm{loc}}	}(\Omega)\equiv\biggl\{
 f\in C^{1}(\Omega):\, f_{|\overline{\Omega'}}\in C^{1,\alpha}(\overline{\Omega'})\ \text{for\ all}\  \Omega'\ \text{such\ that}\ 
 }
\\ \nonumber
&&\qquad\qquad 
 \Omega'\ \text{is\ a\ bounded\ open  subset \ of}\ {\mathbb{R}}^n \,, \overline{\Omega'}\subseteq\Omega  \biggr\}\,,
\end{eqnarray}
 with the family of seminorms
\begin{eqnarray}\label{lem:c1alf2}
\lefteqn{
{\mathcal{P}}_{C^{1,\alpha}_{{\mathrm{loc}}}(\Omega)}\equiv
\biggl\{
\|\cdot\|_{C^{1,\alpha}(\overline{\Omega'})}:\,  
}
\\ \nonumber
&& \qquad\qquad 
 \Omega'\ \text{is\ a\ bounded\ open  \ subset \ of} \ {\mathbb{R}}^n \,, \overline{\Omega'}\subseteq\Omega 
\biggr\}
\end{eqnarray}
is a Fr\'{e}chet space.
\end{lemma}
{\bf Proof.} Since $\Omega$ is the union of an increasing sequence of  sets $\Omega'$ as in (\ref{lem:c1alf2}), the family of seminorms in (\ref{lem:c1alf2}) is equivalent to a countable family of seminorms as in (\ref{lem:c1alf2}) and accordingly,  $C^{1,\alpha}_{	{\mathrm{loc}}	}(\Omega)$ with the above family of seminorms is metrizable. In order to prove that $C^{1,\alpha}_{	{\mathrm{loc}}	}(\Omega)$ is complete it suffices to show that each Cauchy sequence is convergent and to do so 
it suffices to observe that $C^{1,\alpha}(\overline{\Omega'})$ is a Banach space for all $\Omega'$ as in (\ref{lem:c1alf1}).\hfill  $\Box$ 

\vspace{\baselineskip}

Then we can prove the following.
\begin{lemma}\label{lem:c0ademc1a}
  Let   $\alpha\in ]0,1[$. Let $\Omega$ be a bounded open  subset of ${\mathbb{R}}^{n}$. Then $C^{0,\alpha}(\overline{\Omega})_\Delta$ is continuously embedded into the 
  Fr\'{e}chet space $C^{1,\alpha}_{{\mathrm{loc}}}(\Omega)$ with the family of seminorms ${\mathcal{P}}_{C^{1,\alpha}_{{\mathrm{loc}}}(\Omega)}$.
\end{lemma}
{\bf Proof.} By definition of norm in $C^{0,\alpha}(\overline{\Omega})_\Delta$, $C^{0,\alpha}(\overline{\Omega})_\Delta$ is continuously embedded into the Fr\'{e}chet space $C^0(\Omega)$ with the topology of the uniform convergence on the compact subsets of $\Omega$. By the definition of the family 
${\mathcal{P}}_{C^{1,\alpha}_{{\mathrm{loc}}}(\Omega)}$, the 
  Fr\'{e}chet space $C^{1,\alpha}_{{\mathrm{loc}}}(\Omega)$ is also continuously embedded into the Fr\'{e}chet space $C^0(\Omega)$ with the topology of the uniform convergence on the compact subsets of $\Omega$. If $\Omega'$ is as in (\ref{lem:c1alf1}), then there exists an open subset $\Omega_1$ of class $C^\infty$ such that $\overline{\Omega'}\subseteq \Omega_1\subseteq\overline{\Omega_1}\subseteq\Omega$ (cf.~\textit{e.g.}, \cite[Lem.~2.70]{DaLaMu21}).  Then by 
   applying 
  \cite[Lem.~5.7]{La24c} to $\Omega_1$, we conclude that the functions of $C^{0,\alpha}(\overline{\Omega})_\Delta$ are of class $ C^{1,\alpha}_{{\mathrm{loc}}}(\Omega_1)$  and accordingly of class $ C^{1,\alpha}(\overline{\Omega'})$. Hence,  $C^{0,\alpha}(\overline{\Omega})_\Delta\subseteq C^{1,\alpha}_{{\mathrm{loc}}}(\Omega)$. Then we can prove that the embedding of 
$C^{0,\alpha}(\overline{\Omega})_\Delta $ into $C^{1,\alpha}_{{\mathrm{loc}}}(\Omega)$ is continuous by exploiting   the Closed Graph Theorem.\hfill  $\Box$ 

\vspace{\baselineskip}

Next we   introduce a subspace of  $C^{1,\alpha}(\overline{\Omega})$ that we need in the sequel.  
\begin{definition}\label{defn:c1ade}
 Let   $\alpha\in ]0,1[$. Let $\Omega$ be a bounded open  subset of ${\mathbb{R}}^{n}$. Let
 \begin{eqnarray}\label{defn:c1ade1}
C^{1,\alpha}(\overline{\Omega})_\Delta
&\equiv&\biggl\{u\in C^{1,\alpha}(\overline{\Omega}):\,\Delta u\in C^{0,\alpha}(\overline{\Omega})\biggr\}\,,
\\ \nonumber
\|u\|_{ C^{1,\alpha}(\overline{\Omega})_\Delta }
&\equiv& \|u\|_{ C^{1,\alpha}(\overline{\Omega})}
+\|\Delta u\|_{C^{0,\alpha}(\overline{\Omega})}
\qquad\forall u\in C^{1,\alpha}(\overline{\Omega})_\Delta\,.
\end{eqnarray}
\end{definition}
If $\Omega$ is a bounded open subset of ${\mathbb{R}}^{n}$, then 
 $C^{1,\alpha}(\overline{\Omega})$ and $C^{0,\alpha}(\overline{\Omega})$ are Banach spaces and accordingly   $\left(\|u\|_{ C^{1,\alpha}(\overline{\Omega})_\Delta }, \|\cdot \|_{ C^{1,\alpha}(\overline{\Omega})_\Delta }\right)$ is a Banach space. We also note that  if $\Omega$ is a bounded open Lipschitz subset of ${\mathbb{R}}^n$, then 
$C^{1,\alpha}(\overline{\Omega})$ is continuously embedded into $ C^{0,\alpha}(\overline{\Omega})$ and accordingly $C^{1,\alpha}(\overline{\Omega})_\Delta$ is continuously embedded into $ C^{0,\alpha}(\overline{\Omega})_\Delta$. Moreover, by applying \cite[Lem.~5.12]{La24c} to each closed ball  contained in $\Omega$, we deduce that
\begin{equation}\label{eq:c1adc2}
C^{1,\alpha}(\overline{\Omega})_\Delta\subseteq C^2(\Omega)\,.
\end{equation}
 Next we introduce the following classical result on the Green operator for the interior Dirichlet problem. For a proof, we refer for example to \cite[Thm.~4.8]{La24b}.
 \begin{theorem}\label{thm:idwp}
Let $m\in {\mathbb{N}}$, $\alpha\in ]0,1[$. Let $\Omega$ be a bounded open  subset of ${\mathbb{R}}^{n}$ of class $C^{\max\{m,1\},\alpha}$.
Then the map ${\mathcal{G}}_{\Omega,d,+}$  from $C^{m,\alpha}(\partial\Omega)$ to the closed subspace 
 \begin{equation}\label{thm:idwp1}
C^{m,\alpha}_h(\overline{\Omega}) \equiv \{
u\in C^{m,\alpha}(\overline{\Omega}), u\ \text{is\ harmonic\ in}\ \Omega\}
\end{equation}
of $ C^{m,\alpha}(\overline{\Omega})$ that takes $v$ to the only solution $v^\sharp$ of the Dirichlet problem
\begin{equation}\label{defn:cinspo3}
\left\{
\begin{array}{ll}
 \Delta v^\sharp=0 & \text{in}\ \Omega\,,
 \\
v^\sharp_{|\partial\Omega} =v& \text{on}\ \partial\Omega 
\end{array}
\right.
\end{equation}
is a linear homeomorphism.
\end{theorem}
Next we introduce the following two approximation lemmas.
\begin{lemma}\label{lem:apr1a}
 Let $\alpha\in]0,1]$. 
 Let $\Omega$ be a bounded open subset of ${\mathbb{R}}^n$ of class $C^{1,\alpha}$.  If $g\in C^{1,\alpha}(\overline{\Omega})$, then there exists a sequence $\{g_j\}_{j\in {\mathbb{N}}}$ in 
 $C^{\infty}(\overline{\Omega})$ such that
 \begin{equation}\label{lem:apr1a1}
 \sup_{j\in {\mathbb{N}}}\|g_j\|_{
 C^{1,\alpha}(\overline{\Omega}) 
 }<+\infty\,,\qquad
 \lim_{j\to\infty}g_j=g\quad\text{in}\ C^{1,\beta}(\overline{\Omega})\quad \forall\beta\in]0,\alpha[\,.
 \end{equation}
\end{lemma}
For a proof, we refer to \cite[Lem.~A.3]{La24c}. Then we have the following. For a proof, we refer to \cite[Lem.~5.14]{La24c}.

\begin{lemma}\label{lem:apr0ad}
 Let $\alpha\in]0,1]$. 
 Let $\Omega$ be a bounded open subset of ${\mathbb{R}}^n$ of class $C^{1,\alpha}$.  If $u\in C^{0,\alpha}(\overline{\Omega})_\Delta$, then there exists a sequence $\{u_j\}_{j\in {\mathbb{N}}}$ in 
 $C^{1,\alpha}(\overline{\Omega})_\Delta$ such that
 \begin{equation}\label{lem:apr1ad1}
 \sup_{j\in {\mathbb{N}}}\|u_j\|_{
 C^{0,\alpha}(\overline{\Omega})_\Delta
 }<+\infty\,,\qquad
 \lim_{j\to\infty}u_j=u\quad\text{in}\ C^{0,\beta}(\overline{\Omega})_\Delta\quad \forall\beta\in]0,\alpha[\,.
 \end{equation}
\end{lemma}
 Next we plan to introduce the normal derivative of the functions in  $C^{0,\alpha}(\overline{\Omega})_\Delta$ as in \cite{La24c}. To do so,    we introduce the (classical) interior Steklov-Poincar\'{e} operator, \textit{i.e.}, the interior  Dirichlet-to-Neumann map.
\begin{definition}\label{defn:cinspo}
 Let $\alpha\in]0,1[$.  Let  $\Omega$ be a  bounded open subset of ${\mathbb{R}}^{n}$ of class $C^{1,\alpha}$. The classical interior Steklov-Poincar\'{e} operator is defined to be the operator $S_{\Omega,+}$ from
 \begin{equation}\label{defn:cinspo1}
C^{1,\alpha}(\partial\Omega)\quad\text{to}\quad C^{0,\alpha}(\partial\Omega)
\end{equation}
that takes $v\in C^{1,\alpha}(\partial\Omega)$ to the function 
 \begin{equation}\label{defn:cinspo2}
S_{\Omega,+}[v](x)\equiv \frac{\partial  }{\partial\nu}{\mathcal{G}}_{\Omega,d,+}[v](x)\qquad\forall x\in\partial\Omega\,.
\end{equation}
 \end{definition}
   Since   the classical normal derivative is continuous from $C^{1,\alpha}(\overline{\Omega})$ to $C^{0,\alpha}(\partial\Omega)$, the continuity of ${\mathcal{G}}_{\Omega,d,+}$ implies  that $S_{\Omega,+}[\cdot]$ is linear and continuous from 
  $C^{1,\alpha}(\partial\Omega)$ to $C^{0,\alpha}(\partial\Omega)$. Then we have the following definition of \cite[(41)]{La24c}.
  \begin{definition}\label{defn:conoderdedu}
 Let $\alpha\in]0,1[$.  Let  $\Omega$ be a  bounded open subset of ${\mathbb{R}}^{n}$ of class $C^{1,\alpha}$. If  $u\in C^{0}(\overline{\Omega})$ and $\Delta u\in  C^{-1,\alpha}(\overline{\Omega})$, then we define the distributional  normal derivative $\partial_\nu u$
 of $u$ to be the only element of the dual $(C^{1,\alpha}(\partial\Omega))'$ that satisfies the following equality
 \begin{equation}\label{defn:conoderdedu1}
\langle \partial_\nu u ,v\rangle \equiv\int_{\partial\Omega}uS_{\Omega,+}[v]\,d\sigma
+\langle E^\sharp_\Omega[\Delta u],{\mathcal{G}}_{\Omega,d,+}[v]\rangle 
\qquad\forall v\in C^{1,\alpha}(\partial\Omega)\,.
\end{equation}
\end{definition}
The normal derivative of Definition \ref{defn:conoderdedu} extends the classical one in the sense that if $u\in C^{1,\alpha}(\overline{\Omega})$, then under the assumptions  on $\alpha$ and $\Omega$   of Definition \ref{defn:conoderdedu}, we have 
 \begin{equation}\label{lem:conoderdeducl1}
 \langle \partial_\nu u ,v \rangle  \equiv\int_{\partial\Omega}\frac{\partial u}{\partial\nu}v\,d\sigma
  \quad\forall v\in C^{1,\alpha}(\partial\Omega)\,,
\end{equation}
where $\frac{\partial u}{\partial\nu}$ in the right hand side denotes the classical normal derivative of $u$ on $\partial\Omega$ (cf.~\cite[Lem.~5.5]{La24c}). In the sequel, we use the classical symbol $\frac{\partial u}{\partial\nu}$ also for  $\partial_\nu u$ when no ambiguity can arise.\par  
  
Next we introduce the  function space $V^{-1,\alpha}(\partial\Omega)$ on the boundary of $\Omega$ for the normal derivatives of the functions of $C^{0,\alpha}(\overline{\Omega})_\Delta$ as in  \cite[Defn.~13.2, 15.10, Thm.~18.1]{La24b}.
\begin{definition}\label{defn:v-1a}
Let   $\alpha\in ]0,1[$. Let $\Omega$ be a bounded open  subset of ${\mathbb{R}}^{n}$ of class $C^{1,\alpha}$. Let 
\begin{eqnarray}\label{defn:v-1a1}
 \lefteqn{V^{-1,\alpha}(\partial\Omega)\equiv \biggl\{\mu_0+S_{\Omega,+}^t[\mu_1]:\,\mu_0, \mu_1\in C^{0,\alpha}(\partial\Omega)
\biggr\}\,,
}
\\ \nonumber
\lefteqn{
\|\tau\|_{  V^{-1,\alpha}(\partial\Omega) }
\equiv\inf\biggl\{\biggr.
 \|\mu_0\|_{ C^{0,\alpha}(\partial\Omega)  }+\|\mu_1\|_{ C^{0,\alpha}(\partial\Omega)  }
:\,
 \tau=\mu_0+S_{\Omega,+}^t[\mu_1]\biggl.\biggr\}\,,
 }
 \\ \nonumber
 &&\qquad\qquad\qquad\qquad\qquad\qquad\qquad\qquad\qquad
 \forall \tau\in  V^{-1,\alpha}(\partial\Omega)\,,
\end{eqnarray}
where $S_{\Omega,+}^t$ is the transpose map of $S_{\Omega,+}$.
\end{definition}
As shown in \cite[\S 13]{La24b},  $(V^{-1,\alpha}(\partial\Omega), \|\cdot\|_{  V^{-1,\alpha}(\partial\Omega)  })$ is a Banach space. By definition of the norm, $C^{0,\alpha}(\partial\Omega)$ is continuously embedded into $V^{-1,\alpha}(\partial\Omega)$. Moreover, we have the following statement    on the continuity of the normal derivative on $C^{0,\alpha}(\overline{\Omega})_\Delta$  (see also \cite[Prop.~6.6]{La24c}). 
\begin{proposition}\label{prop:ricodnu}
 Let   $\alpha\in ]0,1[$. Let $\Omega$ be a bounded open  subset of 
 ${\mathbb{R}}^{n}$ of class $C^{1,\alpha}$. Then the distributional normal derivative   $\partial_\nu$ is a continuous surjection of $C^{0,\alpha}(\overline{\Omega})_\Delta$ onto $V^{-1,\alpha}(\partial\Omega)$ and there exists $Z\in {\mathcal{L}}\left(V^{-1,\alpha}(\partial\Omega),C^{0,\alpha}(\overline{\Omega})_\Delta\right)$ such that
 \begin{equation}\label{prop:ricodnu1}
\partial_\nu Z[g]=g\qquad\forall g\in V^{-1,\alpha}(\partial\Omega)\,,
\end{equation}
\textit{i.e.}, $Z$ is a right inverse of  $\partial_\nu$.
\end{proposition}

We also mention that the condition of Definition \ref{defn:conoderdedu} admits a different formulation at least in case $u\in C^{0,\alpha}(\overline{\Omega})_\Delta$ (see \cite[Prop.~5.15]{La24c}).
\begin{proposition}\label{prop:node1adeq}
 Let   $\alpha\in ]0,1[$. Let $\Omega$ be a bounded open  subset of 
 ${\mathbb{R}}^{n}$ of class $C^{1,\alpha}$. Let $E_\Omega$ be a linear map from $C^{1,\alpha}(\partial\Omega)$ to $ C^{1,\alpha}(\overline{\Omega})_\Delta$ such that
 \begin{equation}\label{prop:node1adeq0}
 E_\Omega[f]_{|\partial\Omega}=f\qquad\forall f\in C^{1,\alpha}(\partial\Omega)\,.
 \end{equation}
 If $u\in C^{0,\alpha}(\overline{\Omega})_\Delta$, then the distributional  normal derivative $\partial_\nu u$
 of $u$  is characterized by the validity of the following equality
 \begin{eqnarray}\label{prop:node1adeq1}
 \lefteqn{
\langle\partial_\nu u ,v\rangle=\int_{\partial\Omega}u
\frac{\partial}{\partial\nu}E_\Omega[v]
\,d\sigma
}
\\ \nonumber
&&\qquad\qquad
+\langle E^\sharp_\Omega[\Delta u],E_\Omega[v]\rangle -\int_\Omega\Delta (E_\Omega[v]) u\,dx
\qquad\forall v\in C^{1,\alpha}(\partial\Omega)\,.
\end{eqnarray}
\end{proposition}
   We note that in equality (\ref{prop:node1adeq1}) we can use the `extension' operator $E_\Omega$ that we prefer and that accordingly equality (\ref{prop:node1adeq1}) is independent of the specific choice of $E_\Omega$. When we deal with problems for the Laplace operator, a good choice is
 $E_\Omega={\mathcal{G}}_{\Omega,d,+}$, so that the last term in the right hand side of (\ref{prop:node1adeq1}) disappears.  
 In particular, under the assumptions of Proposition \ref{prop:node1adeq}, an extension operator as $E_\Omega$ always exists.
 
 Let $\alpha\in]0,1]$. Let $\Omega$ be a bounded open s
ubset of ${\mathbb{R}}^n$ of class $C^{1,\alpha}$.  We now wish to define the restriction of the distributional derivative $\partial_\nu u $ to a subset $F$ of $\partial\Omega$ that is both open and closed in $\partial\Omega$. To do so, we observe that the linear operator from $C^{1,\alpha}(F)$ to $C^{1,\alpha}(\partial\Omega)$ that is defined by the equality
\[
\stackrel{o}{E}_{F,\partial\Omega}[v]\equiv\left\{
\begin{array}{ll}
 v(x) & \text{if}\ x\in F\,,
 \\
 0 & \text{if}\ x\in (\partial\Omega)\setminus F\,,
\end{array}
\right.
\qquad\forall v\in C^{1,\alpha}(F) 
\]
is linear and continuous and that accordingly the transpose map $\stackrel{o}{E}_{F,\partial\Omega}^t$ is linear and continuous
\[
\text{from}\ 
\left(C^{1,\alpha}(\partial\Omega)\right)'\ 
\text{to}\ \left(C^{1,\alpha}(F)\right)'\,.
\]
Then we are ready to define the restriction of a distributional normal derivative.
\begin{definition}\label{defn:dndf}
 Let   $\alpha\in ]0,1[$. Let $\Omega$ be a bounded open  subset of 
 ${\mathbb{R}}^{n}$ of class $C^{1,\alpha}$. Let $F$ be a subset  of $\partial\Omega$ that is both open and closed in $\partial\Omega$. If $u\in C^{0,\alpha}(\overline{\Omega})_\Delta$, then we say that
 \begin{equation}\label{defn:dndf1}
(\partial_{\nu_\Omega}u)_{|F}\equiv \stackrel{o}{E}_{F,\partial\Omega,}^t\left[
\partial_{\nu_{\Omega}}u
\right] 
\end{equation}
is the restriction to $F$ of the distributional normal derivative $\partial_{\nu_\Omega}u$.
\end{definition}
 Under the assumptions of Definition \ref{defn:dndf},   $(\partial_{\nu_\Omega}u)_{|F}$ is an element of $\left(C^{1,\alpha}(F)\right)'$  and  there exists  a linear (extension) map $E_{\Omega}$ from $C^{1,\alpha}(\partial\Omega )$ to $ C^{1,\alpha}(\overline{\Omega})_\Delta$ such that
 \begin{equation}\label{defn:dndf1a}
 E_{\Omega}[f]_{|\partial\Omega}=f\qquad\forall f\in C^{1,\alpha}(\partial\Omega)\,.
 \end{equation}
 Then the definition of normal derivative on the boundary of 
$\Omega$ implies that
\begin{eqnarray}\label{defn:dndf2}
\lefteqn{
\langle(\partial_{\nu_\Omega}u)_{|F},v\rangle =\langle \partial_{\nu_\Omega}u,\stackrel{o}{E}_{F,\partial\Omega}[v]\rangle
= \int_{F}u\frac{\partial}{\partial\nu_{\Omega}}E_{\Omega}[\stackrel{o}{E}_{F,\partial\Omega}[v]]\,d\sigma
}
\\ \nonumber
&&\qquad\qquad\quad 
+\int_{(\partial\Omega)\setminus F}u\frac{\partial}{\partial\nu_{\Omega}}E_{\Omega}[\stackrel{o}{E}_{F,\partial\Omega}[v]]\,d\sigma
+\langle E^\sharp_{\Omega}[\Delta u],E_{\Omega}[\stackrel{o}{E}_{F,\partial\Omega}[v]]\rangle 
\\ \nonumber
&&\qquad\qquad\quad  
-\int_{\Omega}\Delta (E_{\Omega}[\stackrel{o}{E}_{F,\partial\Omega}[v]]) u\,dx
\qquad\forall v\in C^{1,\alpha}(F)\,,
\end{eqnarray}
(cf.~(\ref{prop:node1adeq1})). Then we also point out the validity of the following elementary remark, that can be readily verified by exploiting the definition of restriction of the normal derivative.
\begin{remark}\label{rem:localization}
 Let   $\alpha\in ]0,1[$. Let $\Omega$ be a bounded open  subset of 
 ${\mathbb{R}}^{n}$ of class $C^{1,\alpha}$. Let $F$ be a subset  of $\partial\Omega$ that is both open and closed in $\partial\Omega$. If $u\in C^{0,\alpha}(\overline{\Omega})_\Delta$, then
 \begin{equation}\label{rem:localization1}
\langle \partial_{\nu_{\Omega}}u,\psi\rangle
=\langle (\partial_{\nu_\Omega}u)_{| F},\psi_{|F}\rangle
+\langle (\partial_{\nu_\Omega}u)_{|(\partial\Omega)\setminus F},\psi_{|(\partial\Omega)\setminus F}\rangle
\qquad\forall \psi\in C^{1,\alpha}(\partial\Omega)\,.
\end{equation}
\end{remark}

Next we assume that $\Omega_1$ and $\Omega_2$ are open subsets of ${\mathbb{R}}^n$ as in Definition \ref{defn:dndf}, that $\Omega_1\subseteq \Omega_2$ and that $F$ is both an open and closed subset of both $\partial\Omega_1$ and  $\partial\Omega_2$ and we ask whether the restriction to $F$ of the normal derivative of a function $u\in C^{0,\alpha}(\overline{\Omega_2})_\Delta$ coincides with the restriction to $F$ of the normal derivative of the restriction of $u$ to $\Omega_1$. In a sense, we ask if the (nonlocal) definition of distributional normal derivative on $F$ would not change if we change the domain away from $F$. We clarify this point by means of the following statement.
\begin{theorem}\label{thm:inrdnd}
 Let   $\alpha\in ]0,1[$. Let $\Omega_1$, $\Omega_2$ be bounded open  subsets of 
 ${\mathbb{R}}^{n}$ of class $C^{1,\alpha}$. Let $\Omega_1\subseteq\Omega_2$.  Let $F$ be both an open and closed subset both of $\partial\Omega_1$ and of $\partial\Omega_2$. Let
 \begin{equation}\label{thm:inrdnd1}
\nu_{\Omega_1}(x)=\nu_{\Omega_2}(x)\qquad\forall x\in F\,.
\end{equation}
Let $F'\equiv (\partial\Omega_1)\setminus F$, $F'\subseteq \Omega_2$. If $u\in C^{0,\alpha}(\overline{\Omega_2})_\Delta$, then 
 \begin{equation}\label{thm:inrdnd2}
(\partial_{\nu_{\Omega_{1}} }u_{|\Omega_1})_{|F}=(\partial_{\nu_{\Omega_{2}} }u)_{|F}
\,.
\end{equation}
\end{theorem}
{\bf Proof.} Since we plan to exploit equality (\ref{defn:dndf2}), we consider the linear extension operator   $E_s$  from the space $C^{1,\alpha}(\partial\Omega_s)$ to $C^{1,\alpha}(\overline{\Omega_s})_\Delta$   for $s\in\{1,2\}$ (cf.~(\ref{defn:dndf1a})). 

 Let  $A_s(u,v)$ denote the right hand side of equality
 (\ref{defn:dndf2}) when   $E_{\Omega}[\cdot]$, $E^\sharp_{\Omega}$, $\stackrel{o}{E}_{F,\partial\Omega}[v]$  are replaced by   $E_s[\cdot]$,  $E^\sharp_{\Omega_s}$,  $\stackrel{o}{E}_{F,\partial\Omega_s}[v]$, respectively, for all $v$ in $C^{1,\alpha}(F)$ and $s\in\{1,2\}$. 
By Lemma \ref{lem:apr0ad},  there exists a sequence $\{u_j\}_{j\in {\mathbb{N}}}$ in 
 $C^{1,\alpha}(\overline{\Omega_2})_\Delta$ such that 
\begin{equation}\label{thm:inrdnd2a}
 \sup_{j\in {\mathbb{N}}}\|u_j\|_{
 C^{0,\alpha}(\overline{\Omega_2})_\Delta
 }<+\infty\,,\ \ 
 \lim_{j\to\infty}u_j=u\quad\text{in}\ C^{0,\beta}(\overline{\Omega_2})_\Delta\ \ \forall\beta\in]0,\alpha[\,.
 \end{equation}
 By the continuity of $\Delta$ from $C^{0,\beta}(\overline{\Omega_s})_\Delta$ to $C^{-1,\beta}(\overline{\Omega_s})$, by the continuity of the extension operator $E_{\Omega_s}^\sharp$ from 
 \[
 C^{-1,\beta}(\overline{\Omega_s})\quad\text{to}\quad\left(C^{1,\beta}(\overline{\Omega_s})\right)'
 \]
 and by the membership of $E_s[\stackrel{o}{E}_{F,\partial\Omega_s}[v]]$ in $C^{1,\alpha}(\overline{\Omega_s})\subseteq C^{1,\beta}(\overline{\Omega_s})$, we have
 \begin{equation}\label{thm:inrdnd3}
 \lim_{j\to\infty}\langle E_{\Omega_s}^\sharp[\Delta u_j],E_s[\stackrel{o}{E}_{F,\partial\Omega_s}[v]]\rangle =\langle E_{\Omega_s}^\sharp[\Delta u],E_s[\stackrel{o}{E}_{F,\partial\Omega_s}[v]]\rangle  \,,
 \end{equation}
 for all $v\in 
 C^{1,\alpha}(F)$, $ s\in\{1,2\}$.   
 
 Then the limiting relations in  (\ref{thm:inrdnd2a}), (\ref{thm:inrdnd3}),  the    H\"{o}lder inequality, equality $E_s[\stackrel{o}{E}_{F,\partial\Omega_s}[v]]_{|(\partial\Omega_s)\setminus F}=0$  for $ s\in\{1,2\}$  and the classical second Green Identity  for the pairs of  functions $u_j$ and $E_2[\stackrel{o}{E}_{F,\partial\Omega_2}[v]]$ in $\Omega_2$
  and then $u_j$ and $E_1[\stackrel{o}{E}_{F,\partial\Omega_1}[v]]$ in $\Omega_1$ (cf.~\textit{e.g.},  \cite[Thm.~4.3]{DaLaMu21} and (\ref{eq:c1adc2})) and assumption (\ref{thm:inrdnd1}) imply that  
 \begin{eqnarray*}
\lefteqn{
\langle (\partial_{\nu_{\Omega_{2}} }u)_{|F}
 ,v\rangle 
}
\\ \nonumber
&& =
A_2(u,v)=\lim_{j\to\infty}A_2(u_j,v)=\lim_{j\to\infty} 
\int_{F}\partial_{\nu_{\Omega_2}}u_jv\,d\sigma=\int_{F}\partial_{\nu_{\Omega_1}}u_jv\,d\sigma
\\ \nonumber
&& 
=\lim_{j\to\infty}A_1(u_j,v)=A_1(u,v)
=\langle (\partial_{\nu_{\Omega_{1}} }u_{|\Omega_1})_{|F},v \rangle 
\ \  \forall v\in 
 C^{1,\alpha}(F) 
\end{eqnarray*}
 and thus the proof is complete.\hfill  $\Box$ 

\vspace{\baselineskip}

We now show that if a function $u$ of class $C^{0,\alpha}(\overline{\Omega})_\Delta$ is classically differentiable around a part $F$ of the boundary as in the Definition \ref{defn:dndf} of restriction, then the restriction of the distributional normal derivative of $u$ to $F$ equals the (distribution that is canonically associated to) the classical normal derivative of $u$ on $F$. 
\begin{proposition}\label{prop:clnore}
 Let   $\alpha\in ]0,1[$. Let $\Omega$ be a bounded open  subset of 
 ${\mathbb{R}}^{n}$ of class $C^{1,\alpha}$. Let $F$ be a subset of $\partial\Omega$ that is both open and closed in $\partial\Omega$. 
 Let $W$ be an open subset of $\Omega$ of class $C^{1,\alpha}$ such that $F$ is a subset of $\partial W$ that is both open and closed in $\partial W$ and
 \[
 \nu_\Omega(x)=\nu_W(x)\qquad\forall x\in F\,.
 \]
 If $u\in C^{0,\alpha}(\overline{\Omega})_\Delta$ and if $u_{|W}\in C^1(\overline{W})$, then 
 $(\partial_{\nu_\Omega}u)_{|F}$ coincides with (the distribution on $F$ that is canonically associated with) the classical normal derivative $\frac{\partial u}{\partial \nu_W}$ on $F$. 
 \end{proposition}
 The validity of Proposition \ref{prop:clnore} follows by applying Theorem \ref{thm:inrdnd} with $\Omega_1=W$ and  $\Omega_2=\Omega$.

\section{A  distributional outward unit normal derivative for the functions of $C^{0,\alpha}_{
{\mathrm{loc}}	}(\overline{\Omega^-})_\Delta$}
\label{sec:doun}
 
 We now plan to define the normal derivative on the boundary for the functions $u$ of $C^{0,\alpha}_{
{\mathrm{loc}}	}(\overline{\Omega^-})_\Delta$ with $\alpha\in]0,1[$ in case $\Omega$ is a bounded open subset of ${\mathbb{R}}^n$ of class $C^{1,\alpha}$. To do so, we  choose $r\in]0,+\infty[$ such that $\overline{\Omega}\subseteq {\mathbb{B}}_n(0,r)$ and we 
plan to define $\partial_{\nu_{\Omega^-}}u$ on $\partial\Omega$ just as the restriction to the
subset $F\equiv\partial\Omega$ of the boundary $(\partial\Omega)\cup(\partial{\mathbb{B}}_n(0,r))$ of ${\mathbb{B}}_n(0,r)\setminus\overline{\Omega}$ of the normal derivative of $u$. In order to shorten our notation, we set
\[
\stackrel{o}{E}_{\partial\Omega,r}[v]\equiv\stackrel{o}{E}_{\partial\Omega, (\partial\Omega)\cup(\partial{\mathbb{B}}_n(0,r))}[v]
\qquad\forall v\in C^{1,\alpha}(\partial\Omega)\,.
\]
 Thus we introduce the following definition.
\begin{definition}\label{defn:endedr}
 Let   $\alpha\in ]0,1[$. Let $\Omega$ be a bounded open  subset of 
 ${\mathbb{R}}^{n}$ of class $C^{1,\alpha}$. Let $r\in]0,+\infty[$ such that $\overline{\Omega}\subseteq {\mathbb{B}}_n(0,r)$. If $u\in C^{0,\alpha}_{
{\mathrm{loc}}	}(\overline{\Omega^-})_\Delta$, then we set
\begin{equation}\label{defn:endedr1}
 \partial_{\nu_{\Omega^-}}u
 \equiv 
 \stackrel{o}{E}_{\partial\Omega,r}^t\left[
\partial_{\nu_{{\mathbb{B}}_n(0,r)\setminus\overline{\Omega}}}u
\right]\,.
\end{equation}
 \end{definition}
 Under the assumptions of Definition \ref{defn:endedr},   $\partial_{\nu_{\Omega^-}}u$ is an element of $\left(C^{1,\alpha}(\partial\Omega)\right)'$.
 Under the assumptions of Definition \ref{defn:endedr}, there exists  a linear (extension) map $E_{{\mathbb{B}}_n(0,r)\setminus\overline{\Omega}}$ from $C^{1,\alpha}((\partial\Omega)\cup(\partial
{\mathbb{B}}_n(0,r)
))$ to $ C^{1,\alpha}(\overline{{\mathbb{B}}_n(0,r)}\setminus\Omega)_\Delta$ such that
 \begin{equation}\label{prop:node1adeq0r}
 E_{{\mathbb{B}}_n(0,r)\setminus\overline{\Omega}}[f]_{|(\partial\Omega)\cup(\partial
{\mathbb{B}}_n(0,r)
)}=f\qquad\forall f\in C^{1,\alpha}((\partial\Omega)\cup(\partial
{\mathbb{B}}_n(0,r)
))\,.
 \end{equation}
Then the definition of normal derivative on the boundary of 
${\mathbb{B}}_n(0,r)\setminus\overline{\Omega}$ implies that
\begin{eqnarray}\label{defn:endedr2}
\lefteqn{
\langle \partial_{\nu_{\Omega^-}}u,v\rangle  
=-\int_{\partial\Omega}u\frac{\partial}{\partial\nu_{\Omega}}E_{{\mathbb{B}}_n(0,r)\setminus\overline{\Omega}}[\stackrel{o}{E}_{\partial\Omega,r}[v]]\,d\sigma
}
\\ \nonumber
&&\qquad\qquad 
+\int_{\partial{\mathbb{B}}_n(0,r)}u\frac{\partial}{\partial\nu_{{\mathbb{B}}_n(0,r)}}E_{{\mathbb{B}}_n(0,r)\setminus\overline{\Omega}}[\stackrel{o}{E}_{\partial\Omega,r}[v]]\,d\sigma
\\ \nonumber
&&\qquad\qquad 
+\langle E^\sharp_{{\mathbb{B}}_n(0,r)\setminus\overline{\Omega}}[\Delta u],E_{{\mathbb{B}}_n(0,r)\setminus\overline{\Omega}}[\stackrel{o}{E}_{\partial\Omega,r}[v]]\rangle 
\\ \nonumber
&&\qquad\qquad 
-\int_{{\mathbb{B}}_n(0,r)\setminus\overline{\Omega}}\Delta (E_{{\mathbb{B}}_n(0,r)\setminus\overline{\Omega}}[\stackrel{o}{E}_{\partial\Omega,r}[v]]) u\,dx
\qquad\forall v\in C^{1,\alpha}(\partial\Omega)\,,
\end{eqnarray}
(cf.~(\ref{prop:node1adeq1})). Now,  we must show that Definition  \ref{defn:endedr} is independent of the choice of $r\in]0,+\infty[$ such that $\overline{\Omega}\subseteq {\mathbb{B}}_n(0,r)$. We do so by means of the following statement.
\begin{proposition}\label{prop:inred}
  Let   $\alpha\in ]0,1[$. Let $\Omega$ be a bounded open  subset of 
 ${\mathbb{R}}^{n}$ of class $C^{1,\alpha}$. Let $u\in C^{0,\alpha}_{
{\mathrm{loc}}	}(\overline{\Omega^-})_\Delta$. Let $r_s\in]0,+\infty[$ be such that $\overline{\Omega}\subseteq {\mathbb{B}}_n(0,r_s)$ for $s\in\{1,2\}$. Then 
 \begin{equation}\label{prop:inred1}
\stackrel{o}{E}_{\partial\Omega,r_1}^t\left[
\partial_{\nu_{{\mathbb{B}}_n(0,r_1)\setminus\overline{\Omega}}}u
\right]=\stackrel{o}{E}_{\partial\Omega,r_2}^t\left[
\partial_{\nu_{{\mathbb{B}}_n(0,r_2)\setminus\overline{\Omega}}}u
\right]\qquad\forall u\in C^{0,\alpha}_{
{\mathrm{loc}}	}(\overline{\Omega^-})_\Delta\,.
\end{equation}
\end{proposition}
{\bf Proof.} There is no loss of generality in assuming that $r_1< r_2$. Then we set
\[
\Omega_s\equiv {\mathbb{B}}_n(0,r_s)\setminus\overline{\Omega}\qquad\forall s\in\{1,2\}\,,
\]
we note that $\partial\Omega$ is both open and closed in the boundary $(\partial\Omega)\cup(\partial{\mathbb{B}}_n(0,r_s))$ of ${\mathbb{B}}_n(0,r_s)\setminus\overline{\Omega}$ for $s\in\{1,2\}$, that
\[
[(\partial\Omega)\cup(\partial{\mathbb{B}}_n(0,r_1))]\setminus(\partial\Omega)
\subseteq \Omega_2\,,
\]
and that accordingly equality (\ref{prop:inred1}) is an immediate consequence of 
Theorem \ref{thm:inrdnd}.\hfill  $\Box$ 

\vspace{\baselineskip}

\begin{remark}\label{rem:conoderdeducl-}
 Let   $\alpha\in ]0,1[$. Let $\Omega$ be a bounded open  subset of 
 ${\mathbb{R}}^{n}$ of class $C^{1,\alpha}$, $u\in C^{1,\alpha}_{
{\mathrm{loc}}	}(\overline{\Omega^-})$. Then    the  
   Definition \ref{defn:endedr} of outward normal derivative and Proposition \ref{prop:clnore} 
   in ${\mathbb{B}}_n(0,r)\setminus\overline{\Omega}$ with $F\equiv\partial\Omega$ imply that
 \begin{equation}\label{rem:conoderdeducl-1}
 \langle \partial_{\nu_{\Omega^-}} u ,v\rangle \equiv\int_{\partial\Omega}\frac{\partial u}{\partial\nu_{\Omega^-}}v\,d\sigma
  \quad\forall v\in C^{1,\alpha}(\partial\Omega)\,,
\end{equation}
where $\frac{\partial u}{\partial\nu_{\Omega^-}}$ in the right hand side denotes the classical  outward normal derivative of $u$ on $\partial\Omega$.  In the sequel, we use the classical symbol $\frac{\partial u}{\partial\nu_{\Omega^-}}$ also for  $\partial_{\nu_{\Omega^-}} u$ when no ambiguity can arise.\par 
\end{remark}

Next we note that if  $u\in C^{0,\alpha}_{{\mathrm{loc}}}(\overline{\Omega^-})$ and $u$ is both harmonic in $\Omega^-$ and harmonic at infinity, then $u\in C^{0,\alpha}_{
{\mathrm{loc}}	}(\overline{\Omega^-})_\Delta$. Accordingly, we now show that the distribution $\partial_{\nu_{\Omega^-}} u $ of Definition \ref{defn:endedr} coincides with the outward normal derivative that has been introduced in \cite[Defn.~6.4]{La24b} for harmonic functions. To do so, we introduce the following classical result on the Green operator for the exterior Dirichlet problem (cf.~\textit{e.g.},  \cite[Thm.~4.11]{La24b}).
\begin{theorem}\label{thm:edwp}
 Let $m\in {\mathbb{N}}$, $\alpha\in ]0,1[$. Let $\Omega$ be a bounded open  subset of ${\mathbb{R}}^{n}$ of class $C^{\max\{m,1\},\alpha}$. The map ${\mathcal{G}}_{\Omega,d,-}$ from $C^{m,\alpha}(\partial\Omega)$ to 
the closed subspace
\begin{equation}\label{thm:edwp1}
C^{m,\alpha}_{bh}(\overline{\Omega^-}) \equiv \{
u\in C^{m,\alpha}_b(\overline{\Omega^-}), u\ \text{is\ harmonic\ in}\ \Omega^- 
 \text{and\ is\ harmonic\ at}\ \infty\}
\end{equation}
of $C^{m,\alpha}_{b}(\overline{\Omega^-})$ that takes $v$ to the only solution $v^\sharp$ of   in $ C^{m,\alpha}_{b}(\overline{\Omega^-})$ of the exterior Dirichlet problem
\begin{equation}\label{defn:censpo3}
\left\{
\begin{array}{ll}
 \Delta v^\sharp=0 & \text{in}\ \Omega^-\,,
 \\
v^\sharp_{|\partial\Omega} =v& \text{on}\ \partial\Omega\,,
\\
v^\sharp\  \text{is\ harmonic\ at}\ \infty   &
\end{array}
\right.
\end{equation}
is a linear homeomorphism.  
\end{theorem} 
We are now ready to prove  the following statement. 
\begin{proposition}\label{prop:coenoder}
  Let   $\alpha\in ]0,1[$. Let $\Omega$ be a bounded open  subset of 
 ${\mathbb{R}}^{n}$ of class $C^{1,\alpha}$. Let $r\in]0,+\infty[$ be such that $\overline{\Omega}\subseteq {\mathbb{B}}_n(0,r)$.   If $u\in C^{0,\alpha}_{{\mathrm{loc}}}(\overline{\Omega^-})$ and $u$ is both harmonic in $\Omega^-$ and harmonic at infinity, then
\begin{equation}\label{prop:coenoder1}
\int_{\partial\Omega}u\frac{\partial}{\partial\nu_{\Omega^-}} {\mathcal{G}}_{\Omega,d,-}[v]\,d\sigma
=\langle \stackrel{o}{E}_{\partial\Omega,r}^t\left[
\partial_{\nu_{{\mathbb{B}}_n(0,r)\setminus\overline{\Omega}}}u
\right],v\rangle \qquad\forall v\in C^{1,\alpha}(\partial\Omega)\,,
\end{equation}
\textit{i.e.},  the definition of $\partial_{\nu_{\Omega^-}}u$ in (\ref{defn:endedr1}) (or its equivalent form (\ref{defn:endedr2})) coincides with that of \cite[Defn.~6.4]{La24b}. 
 \end{proposition}
{\bf Proof.} Let $E$ be a linear extension operator  from $C^{1,\alpha}(\partial\Omega)$ to $C^{1,\alpha}(\overline{{\mathbb{B}}_n(0,r)}\setminus\Omega)_\Delta$. A known approximation property implies that there exists a sequence $\{f_j\}_{j\in {\mathbb{N}}}$ in $C^{1,\alpha}(\partial\Omega)$   that converges to $u_{|\partial\Omega}$ in the $ C^{0,\beta}(\partial\Omega)$-norm for all $\beta\in]0,\alpha[$ and that is bounded in the $ C^{0,\alpha}(\partial\Omega)$-norm (cf.~\textit{e.g.}, \cite[Lem.~A.25]{La24b}). Next we set
\[
u_j\equiv {\mathcal{G}}_{\Omega,d,-}[f_j]\qquad\forall j\in {\mathbb{N}}\,.
\]
By the classical Theorem \ref{thm:edwp} on the Green operator, we have
\begin{eqnarray}\label{prop:coenoder1a}
&&u_j\in C^{1,\alpha}_{bh}(\overline{\Omega^-}) \qquad\forall j\in {\mathbb{N}}\,,
\quad
\sup_{j\in {\mathbb{N}}}\|u_j\|_{C^{0,\alpha}_{bh}(\overline{\Omega^-})}
<+\infty\,,
\\ \nonumber
&&\lim_{j\to\infty}u_j=u
\ \text{in}\ C^{0,\beta}_{bh}(\overline{\Omega^-})\ \qquad\forall \beta\in]0,\alpha[\,.
\end{eqnarray}
 Clearly
 \begin{eqnarray}\label{prop:coenoder1b}
\lefteqn{
\langle E^\sharp_{{\mathbb{B}}_n(0,r)\setminus\overline{\Omega}}[\Delta u_j],E_{{\mathbb{B}}_n(0,r)\setminus\overline{\Omega}}[\stackrel{o}{E}_{\partial\Omega,r}[v]]\rangle  
}
\\ \nonumber
&&\qquad 
=0 =\langle E^\sharp_{{\mathbb{B}}_n(0,r)\setminus\overline{\Omega}}[\Delta u],E_{{\mathbb{B}}_n(0,r)\setminus\overline{\Omega}}[\stackrel{o}{E}_{\partial\Omega,r}[v]]\rangle  
\qquad\forall v\in 
 C^{1,\alpha}(\partial\Omega)\,.
\end{eqnarray}
Now let $A (u,v)$ denote the right hand side of 
(\ref{defn:endedr2}).  Then the limiting relations in (\ref{prop:coenoder1a}), the    H\"{o}lder inequality,
the equalities in (\ref{prop:coenoder1b}),    equality 
\[
E_{{\mathbb{B}}_n(0,r)\setminus\overline{\Omega}}[\stackrel{o}{E}_{\partial\Omega,r}[v]]_{|\partial{\mathbb{B}}_{n}(0,r)}=0
\] and the classical second Green Identity  in $ {\mathbb{B}}_n(0,r)\setminus\overline{\Omega}$   for the pairs of  functions $u_j$ and $E_{{\mathbb{B}}_n(0,r)\setminus\overline{\Omega}}[\stackrel{o}{E}_{\partial\Omega,r}[v]]$ (cf.~\textit{e.g.},  \cite[Cor.~4.3]{DaLaMu21} and (\ref{eq:c1adc2})) imply that
  \begin{equation}\label{prop:coenoder2}
A(u,v)=\lim_{j\to\infty}A(u_j,v)=\lim_{j\to\infty} 
\int_{\partial\Omega}\partial_{\nu_{\Omega^-}}u_jv\,d\sigma\,.
\end{equation}
Then the second Green Identity for harmonic functions in $\Omega^-$ that are harmonic at infinity    for the pairs of  functions $u_j$ and ${\mathcal{G}}_{\Omega,d,-}[v]$ implies that
 \begin{eqnarray}\label{prop:coenoder2a}
\lefteqn{\lim_{j\to\infty} 
\int_{\partial\Omega}\partial_{\nu_{\Omega^-}}u_jv\,d\sigma 
=\lim_{j\to\infty} 
\int_{\partial\Omega}\partial_{\nu_{\Omega^-}}u_j{\mathcal{G}}_{\Omega,d,-}[v]\,d\sigma 
}
\\ \nonumber
&&\qquad
=
\lim_{j\to\infty} 
 \int_{\partial\Omega}u_j\partial_{\nu_{\Omega^-}}{\mathcal{G}}_{\Omega,d,-}[v]\,d\sigma 
 = \int_{\partial\Omega}u\partial_{\nu_{\Omega^-}}{\mathcal{G}}_{\Omega,d,-}[v]\,d\sigma 
\end{eqnarray}
(cf.~\textit{e.g.},  \cite[Cor.~4.8]{DaLaMu21}). Hence, (\ref{prop:coenoder1}) holds true and the proof is complete.\hfill  $\Box$ 

\vspace{\baselineskip}

In order to state a characterization of the space $V^{-1,\alpha}(\partial\Omega)$, we introduce the following definition.
\begin{definition}\label{defn:censpo}
 Let $\alpha\in]0,1[$.  Let  $\Omega$ be a  bounded open subset of ${\mathbb{R}}^{n}$ of class $C^{1,\alpha}$. The classical interior Steklov-Poincar\'{e} operator is defined to be the operator $S_-$ from
 \begin{equation}\label{defn:censpo1}
C^{1,\alpha}(\partial\Omega)\quad\text{to}\quad C^{0,\alpha}(\partial\Omega)
\end{equation}
 takes $v\in C^{1,\alpha}(\partial\Omega)$ to the function 
 \begin{equation}\label{defn:censpo2}
S_{\Omega,-}[v](x)\equiv \frac{\partial}{\partial\nu_{\Omega^-}}{\mathcal{G}}_{\Omega,d,-}[v](x)\qquad\forall x\in\partial\Omega\,,
\end{equation}
where $\nu_{\Omega^-}=-\nu_\Omega$ is the outward unit normal to $\Omega^-$ on $\partial\Omega=\partial\Omega^-$.
 \end{definition}
 Since   the classical normal derivative is continuous from $C^{1,\alpha}_b(\overline{\Omega^-})$ to $C^{0,\alpha}(\partial\Omega)$, the continuity of ${\mathcal{G}}_{\Omega,d,-}$ implies  that the operator $S_{\Omega,-}[\cdot]$ is linear and continuous from 
  $C^{1,\alpha}(\partial\Omega)$ to $C^{0,\alpha}(\partial\Omega)$.  We are now ready to state the following  characterization of the space $V^{-1,\alpha}(\partial\Omega)$ (see \cite[Defn.~13.2, 15.10, Thm.~18.1]{La24b}).
\begin{theorem}\label{thm:v-1a-}
 Let   $\alpha\in ]0,1[$. Let $\Omega$ be a bounded open  subset of ${\mathbb{R}}^{n}$ of class $C^{1,\alpha}$. Then
 \begin{equation}\label{thm:v-1a-1}
 V^{-1,\alpha}(\partial\Omega)=\biggl\{\mu_0+S_{\Omega,-}^t[\mu_1]:\,\mu_0, \mu_1\in C^{0,\alpha}(\partial\Omega)
\biggr\}\,,
\end{equation}
Moreover, the norm
\begin{eqnarray}\label{thm:v-1a-2}
\lefteqn{
\|\tau\|_{  V^{-1,\alpha}(\partial\Omega) -}
\equiv\inf\biggl\{\biggr.
 \|\mu_0\|_{ C^{0,\alpha}(\partial\Omega)  }+\|\mu_1\|_{ C^{0,\alpha}(\partial\Omega)  }
:\,
 }
 \\ \nonumber
 &&\qquad\qquad\qquad\qquad\qquad\qquad\ \tau=\mu_0+S_{\Omega,-}^t[\mu_1]\biggl.\biggr\}
\qquad 
 \forall \tau\in  V^{-1,\alpha}(\partial\Omega)\,,
\end{eqnarray}
where $S_{\Omega,-}^t$ is the transpose map of $S_{\Omega,-}$ is equivalent to the norm of $V^{-1,\alpha}(\partial\Omega)$ that has been defined by  (\ref{defn:v-1a1}).
\end{theorem}

Next we want to show that the distributional normal derivative $\partial_{\nu_{\Omega^-}}$ of
Definition \ref{defn:endedr} is linear and continuous from  $C^{0,\alpha}_{
{\mathrm{loc}}	}(\overline{\Omega^-})_\Delta$ to the space $V^{-1,\alpha}(\partial\Omega)$. By Propostion \ref{prop:ricodnu} and formula (\ref{defn:endedr1}) it suffices to show that $\stackrel{o}{E}_{\partial\Omega,r}^t$ is linear and continuous from  $V^{-1,\alpha}((\partial\Omega)\cup (\partial{\mathbb{B}}_n(0,r)))$ to $V^{-1,\alpha}(\partial\Omega)$. We do so by means of the following statement.

\begin{proposition}\label{prop:reco-1a}
  Let   $\alpha\in ]0,1[$. Let $\Omega$ be a bounded open  subset of 
 ${\mathbb{R}}^{n}$ of class $C^{1,\alpha}$. Let $r\in]0,+\infty[$ such that $\overline{\Omega}\subseteq {\mathbb{B}}_n(0,r)$. Then the map $\stackrel{o}{E}_{\partial\Omega,r}^t$ is a linear and continuous surjection from  $V^{-1,\alpha}((\partial\Omega)\cup (\partial{\mathbb{B}}_n(0,r)))$ to $V^{-1,\alpha}(\partial\Omega)$. 
 
 Moreover the transpose map $r_{|\partial\Omega}^t$ of the restriction operator $r_{|\partial\Omega}$ from the space 
 $C^{1,\alpha}((\partial\Omega)\cup (\partial{\mathbb{B}}_n(0,r)))$ to $C^{1,\alpha}(\partial\Omega)$ is a linear and continuous injection from $V^{-1,\alpha}(\partial\Omega)$ to $V^{-1,\alpha}((\partial\Omega)\cup (\partial{\mathbb{B}}_n(0,r)))$ and
 \begin{equation}\label{prop:reco-1a1}
 \stackrel{o}{E}_{\partial\Omega,r}^t\circ r_{|\partial\Omega}^t[\tau]=\tau\qquad\forall \tau\in V^{-1,\alpha}(\partial\Omega)\,,
\end{equation}
\textit{i.e.}, $r_{|\partial\Omega}^t$ is a right inverse of  $\stackrel{o}{E}_{\partial\Omega,r}^t$.
 \end{proposition}
{\bf Proof.} We first prove that if $\tau\in  V^{-1,\alpha}((\partial\Omega)\cup (\partial{\mathbb{B}}_n(0,r)))$, then 
$\stackrel{o}{E}_{\partial\Omega,r}^t[\tau]$ belongs to $ V^{-1,\alpha}(\partial\Omega)$. By Theorem \ref{thm:v-1a-}, there exist $\mu_0$, $\mu_1\in C^{0,\alpha}((\partial\Omega)\cup (\partial{\mathbb{B}}_n(0,r)))$ such that
 \begin{equation}\label{prop:reco-1a2}
 \tau=\mu_0+S_{{\mathbb{B}}_n(0,r)\setminus\overline{\Omega},-}^t[\mu_1]\,.
\end{equation}
By the definition of $\stackrel{o}{E}_{\partial\Omega,r}$, we have
 \begin{equation}\label{prop:reco-1a3}
  \stackrel{o}{E}_{\partial\Omega,r}^t[\mu_0]=\mu_{0|\partial\Omega}\in 
  C^{0,\alpha}(\partial\Omega)\subseteq 
  V^{-1,\alpha}(\partial\Omega)\,.
\end{equation}
We now turn to show that
\begin{equation}\label{prop:reco-1a3a}
 \stackrel{o}{E}_{\partial\Omega,r}^tS_{{\mathbb{B}}_n(0,r)\setminus\overline{\Omega},-}^t[\mu_1]
 =-S_{\Omega,+}^t[\mu_{1|\partial\Omega}]\qquad\text{on}\ \partial\Omega\,.
\end{equation}
Let $\psi\in C^{1,\alpha}(\partial\Omega)$. Then
\begin{eqnarray}\label{prop:reco-1a4}
\lefteqn{
\langle -S_{\Omega,+}^t[\mu_{1|\partial\Omega}],\psi\rangle 
}
\\ \nonumber
&&\qquad
=-\langle \mu_{1|\partial\Omega},S_{\Omega,+}[\psi]\rangle 
=-\langle \mu_{1|\partial\Omega},\frac{\partial}{\partial\nu_\Omega}{\mathcal{G}}_{\Omega,d,+}[\psi]\rangle \,.
\end{eqnarray}
By the definition of the Green operator, we have
 \[
 {\mathcal{G}}_{\Omega,d,+}[\psi]
 =
  {\mathcal{G}}_{
  {\mathbb{B}}_n(0,r)\setminus\overline{\Omega},d,-
  }[\stackrel{o}{E}_{\partial\Omega,r}[\psi]]_{|\overline{\Omega}}
  \]
and accordingly,
\begin{eqnarray} \nonumber
&&
\frac{\partial}{\partial\nu_\Omega}{\mathcal{G}}_{\Omega,d,+}[\psi](x)
=-\frac{\partial}{\partial\nu_{{\mathbb{B}}_n(0,r)\setminus\overline{\Omega}}}{\mathcal{G}}_{{\mathbb{B}}_n(0,r)\setminus\overline{\Omega},d,-}[\stackrel{o}{E}_{\partial\Omega,r}[\psi]](x)
\quad\forall x\in\partial\Omega\,,
\\ \label{prop:reco-1a5}
&& 
S_{\Omega,+}[\psi](x)=-S_{
{\mathbb{B}}_n(0,r)\setminus\overline{\Omega},-}
[\stackrel{o}{E}_{\partial\Omega,r}[\psi]](x)
 \qquad\forall x\in\partial\Omega\,.
\end{eqnarray}
Hence, equality (\ref{prop:reco-1a4}) implies that
\begin{equation}\label{prop:reco-1a6}
\langle -S_{\Omega,+}^t[\mu_{1|\partial\Omega}],\psi\rangle 
=\langle \mu_{1|\partial\Omega},S_{
{\mathbb{B}}_n(0,r)\setminus\overline{\Omega},-}
[\stackrel{o}{E}_{\partial\Omega,r}[\psi]]_{|\partial\Omega}\rangle \,.
\end{equation}
Next we note that 
\[
{\mathcal{G}}_{{\mathbb{B}}_n(0,r)\setminus\overline{\Omega},d,-}[\stackrel{o}{E}_{\partial\Omega,r}[\psi]](x)=0
\qquad\forall x\in {\mathbb{R}}^n\setminus{\mathbb{B}}_n(0,r)
\]
and that accordingly
\[
\frac{\partial}{\partial\nu_{{\mathbb{B}}_n(0,r)\setminus\overline{\Omega}}}{\mathcal{G}}_{{\mathbb{B}}_n(0,r)\setminus\overline{\Omega},d,-}[\stackrel{o}{E}_{\partial\Omega,r}[\psi]](x)=0\qquad \forall x\in \partial{\mathbb{B}}_n(0,r)\,.
\]
Then equality (\ref{prop:reco-1a4}) implies that
\begin{eqnarray}\label{prop:reco-1a6a}
\lefteqn{
\langle -S_{\Omega,+}^t[\mu_{1|\partial\Omega}],\psi\rangle 
}
\\ \nonumber
&&=\int_{
(\partial\Omega)\cup (\partial{\mathbb{B}}_n(0,r))
}\mu_1\frac{\partial}{\partial\nu_{{\mathbb{B}}_n(0,r)\setminus\overline{\Omega}}}{\mathcal{G}}_{{\mathbb{B}}_n(0,r)\setminus\overline{\Omega},d,-}[\stackrel{o}{E}_{\partial\Omega,r}[\psi]]\,d\sigma
\\ \nonumber
&&
=\langle \mu_1,S_{{\mathbb{B}}_n(0,r)\setminus\overline{\Omega},-}[\stackrel{o}{E}_{\partial\Omega,r}[\psi]]\rangle 
=\langle \stackrel{o}{E}_{\partial\Omega,r}^tS_{{\mathbb{B}}_n(0,r)\setminus\overline{\Omega},-}^t[\mu_1],\psi\rangle \,.
\end{eqnarray}
Hence, equality (\ref{prop:reco-1a3a}) holds true. Since $S_{\Omega,+}^t[\mu_{1|\partial\Omega}]\in V^{-1,\alpha}(\partial\Omega)$, we conclude that $\stackrel{o}{E}_{\partial\Omega,r}^tS_{{\mathbb{B}}_n(0,r)\setminus\overline{\Omega},-}^t[\mu_1]\in V^{-1,\alpha}(\partial\Omega)$ and thus the membership of (\ref{prop:reco-1a3}) implies that $\stackrel{o}{E}_{\partial\Omega,r}^t[\tau]\in V^{-1,\alpha}(\partial\Omega)$.

Next we show the continuity statement of $\stackrel{o}{E}_{\partial\Omega,r}^t$. Since $\stackrel{o}{E}_{\partial\Omega,r}$ is linear and continuous from $C^{1,\alpha}(\partial\Omega)$ to $C^{1,\alpha}((\partial\Omega)\cup (\partial{\mathbb{B}}_n(0,r)))$, we have 
\[
\stackrel{o}{E}_{\partial\Omega,r}^t\in {\mathcal{L}}\left(
\left(C^{1,\alpha}((\partial\Omega)\cup (\partial{\mathbb{B}}_n(0,r)))\right)',
\left(C^{1,\alpha}(\partial\Omega)\right)'
\right)\,.
\]
Then the continuity of the embedding of  $V^{-1,\alpha}((\partial\Omega)\cup (\partial{\mathbb{B}}_n(0,r)))$  into the space 
$\left(C^{1,\alpha}((\partial\Omega)\cup (\partial{\mathbb{B}}_n(0,r)))\right)'$ implies that
\begin{equation}\label{prop:reco-1a7}
\stackrel{o}{E}_{\partial\Omega,r}^t\in {\mathcal{L}}\left(
V^{-1,\alpha}((\partial\Omega)\cup (\partial{\mathbb{B}}_n(0,r))),
\left(C^{1,\alpha}(\partial\Omega)\right)'
\right)\,.
\end{equation}
Since we have proved that
\[
\stackrel{o}{E}_{\partial\Omega,r}^t[V^{-1,\alpha}((\partial\Omega)\cup (\partial{\mathbb{B}}_n(0,r)))]\subseteq V^{-1,\alpha}(\partial\Omega)
\]
and  $V^{-1,\alpha}(\partial\Omega)$ is continuously embedded into $\left(C^{1,\alpha}(\partial\Omega)\right)'$, the Closed Graph Theorem implies that
\begin{equation}\label{prop:reco-1a8}
\stackrel{o}{E}_{\partial\Omega,r}^t\in {\mathcal{L}}\left(
V^{-1,\alpha}((\partial\Omega)\cup (\partial{\mathbb{B}}_n(0,r))),
V^{-1,\alpha}(\partial\Omega)
\right)\,.
\end{equation}
Next we prove that  if $\varsigma\in  V^{-1,\alpha}(\partial\Omega)$, then 
$r_{|\partial\Omega}^t[\varsigma]\in  V^{-1,\alpha}((\partial\Omega)\cup (\partial{\mathbb{B}}_n(0,r)))$. By the definition of 
$V^{-1,\alpha}(\partial\Omega)$, there exist $f_0$, $f_1\in C^{0,\alpha}(\partial\Omega )$ such that
\[
\varsigma=f_0+S_{\Omega,+}^t[f_1]\,.
\]
By the definition of $\stackrel{o}{E}_{\partial\Omega,r}$, we have
\begin{eqnarray}\label{prop:reco-1a9}
\lefteqn{
\langle r_{|\partial\Omega}^t[f_0],\varphi \rangle =\langle f_0,r_{|\partial\Omega}[\varphi]\rangle 
}
\\
\nonumber
&&\qquad=
\int_{\partial\Omega}f_0r_{|\partial\Omega}[\varphi]\,d\sigma
=\int_{(\partial\Omega)\cup (\partial{\mathbb{B}}_n(0,r))}\stackrel{o}{E}_{\partial\Omega,r}[f_0]\varphi\,d\sigma
=\langle \stackrel{o}{E}_{\partial\Omega,r}[f_0],\varphi\rangle 
\end{eqnarray}
for all $\varphi\in C^{1,\alpha}((\partial\Omega)\cup (\partial{\mathbb{B}}_n(0,r)))$. We now turn to compute
$r_{|\partial\Omega}^tS_{\Omega,+}^t[f_1]$. Let $\varphi\in C^{1,\alpha}((\partial\Omega)\cup (\partial{\mathbb{B}}_n(0,r)))$. Then
\begin{eqnarray}\label{prop:reco-1a10}
\lefteqn{
\langle r_{|\partial\Omega}^tS_{\Omega,+}^t[f_1],\varphi\rangle 
}
\\ \nonumber
&&\qquad
=\langle f_1,S_{\Omega,+}[r_{|\partial\Omega}\varphi]\rangle =\int_{\partial\Omega}f_1\frac{\partial}{\partial\nu_\Omega}{\mathcal{G}}_{\Omega,d,+}[r_{|\partial\Omega}\varphi]\,d\sigma\,.
\end{eqnarray}
 By the definition of the Green operator, we have
 \[
 {\mathcal{G}}_{\Omega,d,+}[r_{|\partial\Omega}\varphi](x)
 =
  {\mathcal{G}}_{
  {\mathbb{B}}_n(0,r)\setminus\overline{\Omega},d,-
  }[\varphi](x)\qquad\forall x\in\overline{\Omega}
  \]
and accordingly,
\begin{eqnarray} \nonumber
&&
\frac{\partial}{\partial\nu_\Omega}{\mathcal{G}}_{\Omega,d,+}[r_{|\partial\Omega}\varphi](x)
=
-\frac{\partial}{\partial\nu_{{\mathbb{B}}_n(0,r)\setminus\overline{\Omega}}}{\mathcal{G}}_{{\mathbb{B}}_n(0,r)\setminus\overline{\Omega},d,-}[\varphi](x)
\qquad\forall x\in\partial\Omega\,,
\\ \label{prop:reco-1a11}
&& 
S_{\Omega,+}[r_{|\partial\Omega}\varphi](x)=-S_{
{\mathbb{B}}_n(0,r)\setminus\overline{\Omega},-}
[\varphi](x)
 \qquad\forall x\in\partial\Omega\,.
\end{eqnarray}
Then equality (\ref{prop:reco-1a10}) implies that
\begin{eqnarray}\label{prop:reco-1a12}
\lefteqn{
\langle r_{|\partial\Omega}^tS_{\Omega,+}^t[f_1],\varphi\rangle 
=-\int_{\partial\Omega}f_1S_{
{\mathbb{B}}_n(0,r)\setminus\overline{\Omega},-}
[\varphi] \,d\sigma
}
\\ \nonumber
&&\qquad\qquad\qquad\qquad 
=-\int_{(\partial\Omega)\cup (\partial{\mathbb{B}}_n(0,r))}\stackrel{o}{E}_{\partial\Omega,r}[f_1]S_{
{\mathbb{B}}_n(0,r)\setminus\overline{\Omega},-}
[\varphi] \,d\sigma
\\ \nonumber
&&\qquad\qquad\qquad\qquad
=-\langle S_{
{\mathbb{B}}_n(0,r)\setminus\overline{\Omega},-}^t[\stackrel{o}{E}_{\partial\Omega,r}[f_1]],\varphi\rangle \,.
\end{eqnarray}
Then equalities (\ref{prop:reco-1a9}) and (\ref{prop:reco-1a12}) imply that
\begin{equation}\label{prop:reco-1a13}
r_{|\partial\Omega}^t[\varsigma]=r_{|\partial\Omega}^t[f_0]+ 
r_{|\partial\Omega}^tS_{\Omega,+}^t[f_1]=\stackrel{o}{E}_{\partial\Omega,r}[f_0]-S_{
{\mathbb{B}}_n(0,r)\setminus\overline{\Omega},-}^t[\stackrel{o}{E}_{\partial\Omega,r}[f_1]]
\end{equation}
Since $\stackrel{o}{E}_{\partial\Omega,r}[f_0]$, $-\stackrel{o}{E}_{\partial\Omega,r}[f_1]\in C^{0,\alpha}((\partial\Omega)\cup (\partial{\mathbb{B}}_n(0,r)))$, the characterization Theorem \ref{thm:v-1a-} 
  implies that $r_{|\partial\Omega}^t[\varsigma]\in V^{-1,\alpha}((\partial\Omega)\cup (\partial{\mathbb{B}}_n(0,r)))$.

Next we show the continuity statement of $r_{|\partial\Omega}^t$. Since $r_{|\partial\Omega}$ is linear and continuous from $C^{1,\alpha}((\partial\Omega)\cup (\partial{\mathbb{B}}_n(0,r)))$ to $C^{1,\alpha}(\partial\Omega)$, we have 
\[
r_{|\partial\Omega}^t\in {\mathcal{L}}\left(
 \left(C^{1,\alpha}(\partial\Omega)\right)',
\left(C^{1,\alpha}((\partial\Omega)\cup (\partial{\mathbb{B}}_n(0,r)))\right)'
\right)\,.
\]
Then the continuity of the embedding of  $V^{-1,\alpha}(\partial\Omega))$  into
$\left(C^{1,\alpha}( \partial\Omega)\right)'$ implies that
\begin{equation}\label{prop:reco-1a14}
r_{|\partial\Omega}^t\in {\mathcal{L}}\left(
V^{-1,\alpha}(\partial\Omega),
\left(C^{1,\alpha}((\partial\Omega)\cup (\partial{\mathbb{B}}_n(0,r)))\right)'
\right)\,.
\end{equation}
Since we have proved that
\[
r_{|\partial\Omega}^t[V^{-1,\alpha}(\partial\Omega)]\subseteq  V^{-1,\alpha}((\partial\Omega)\cup (\partial{\mathbb{B}}_n(0,r)))
\]
and  the space $V^{-1,\alpha}((\partial\Omega)\cup (\partial{\mathbb{B}}_n(0,r)))$ is continuously embedded into the space  $\left(C^{1,\alpha}((\partial\Omega)\cup (\partial{\mathbb{B}}_n(0,r))))\right)'$, the Closed Graph Theorem implies that
\begin{equation}\label{prop:reco-1a15}
r_{|\partial\Omega}^t\in {\mathcal{L}}\left(V^{-1,\alpha}(\partial\Omega),
V^{-1,\alpha}((\partial\Omega)\cup (\partial{\mathbb{B}}_n(0,r))) 
 \right)\,.
\end{equation}
Since $r_{|\partial\Omega}\circ \stackrel{o}{E}_{\partial\Omega,r}$ is the identity map in $C^{1,\alpha}(\partial\Omega)$, by taking the transpose maps we conclude that equality (\ref{prop:reco-1a1}) holds true. Then equality (\ref{prop:reco-1a1}) implies that the map $\stackrel{o}{E}_{\partial\Omega,r}^t$  from $V^{-1,\alpha}((\partial\Omega)\cup (\partial{\mathbb{B}}_n(0,r)))$ to $V^{-1,\alpha}(\partial\Omega)$ is a surjection
and that the map $r_{|\partial\Omega}^t$ from $V^{-1,\alpha}(\partial\Omega)$ to $V^{-1,\alpha}((\partial\Omega)\cup (\partial{\mathbb{B}}_n(0,r)))$  is injective.\hfill  $\Box$ 

\vspace{\baselineskip}

 We are now ready to prove the following theorem.
\begin{proposition}\label{prop:recodnu}
 Let   $\alpha\in ]0,1[$. Let $\Omega$ be a bounded open  subset of 
 ${\mathbb{R}}^{n}$ of class $C^{1,\alpha}$.  Then the distributional normal derivative   $\partial_{\nu_{\Omega^-}}$ is a continuous surjection from $C^{0,\alpha}_{
{\mathrm{loc}}	}(\overline{\Omega^-})_\Delta$ onto $V^{-1,\alpha}(\partial\Omega)$ and there exists a linear and continuous  operator   $Z_-$ from  $ V^{-1,\alpha}(\partial\Omega)$ to $C^{0,\alpha}_{
{\mathrm{loc}}	}(\overline{\Omega^-})_\Delta $ such that
 \begin{equation}\label{prop:recodnu1}
\partial_{\nu_{\Omega^-}} Z_-[g]=g\qquad\forall g\in V^{-1,\alpha}(\partial\Omega)\,,
\end{equation}
\textit{i.e.}, $Z_-$ is a right inverse of  $\partial_{\nu_{\Omega^-}}$. (See Lemma \ref{lem:c1alcof} for the topology of $C^{0,\alpha}_{
{\mathrm{loc}}	}(\overline{\Omega^-})_\Delta $).
\end{proposition}
{\bf Proof.} Let $r\in]0,+\infty[$ be such that $\overline{\Omega}\subseteq {\mathbb{B}}_n(0,r)$. By the Definition \ref{defn:endedr} of distrubutional outward normal derivative,  
by the continuity of the restriction map from $C^{0,\alpha}_{
{\mathrm{loc}}	}(\overline{\Omega^-})_\Delta$ to $C^{0,\alpha}(\overline{{\mathbb{B}}_n(0,r)}\setminus\Omega)_\Delta$, 
by the continuity of 
$\partial_{\nu_{{\mathbb{B}}_n(0,r)\setminus\overline{\Omega}}}$ of Proposition \ref{prop:ricodnu} 
from $C^{0,\alpha}(\overline{{\mathbb{B}}_n(0,r)}\setminus\Omega)_\Delta$ to $V^{-1,\alpha}((\partial\Omega)\cup (\partial{\mathbb{B}}_n(0,r)))$ and
 by the continuity of $ \stackrel{o}{E}_{\partial\Omega,r}^t$ of Proposition \ref{prop:reco-1a}, we conclude that $\partial_{\nu_{\Omega^-}}$ is linear and continuous from 
$C^{0,\alpha}_{
{\mathrm{loc}}	}(\overline{\Omega^-})_\Delta$ to  $V^{-1,\alpha}(\partial\Omega)$.

Next, we turn to show the existence of $Z_-$. By Proposition \ref{prop:ricodnu} and Proposition \ref{prop:reco-1a}, there exists
$Z\in {\mathcal{L}}\left(V^{-1,\alpha}((\partial\Omega)\cup (\partial{\mathbb{B}}_n(0,r))),C^{0,\alpha}(\overline{{\mathbb{B}}_n(0,r)}\setminus\Omega)_\Delta\right)$ such that
 \begin{equation}\label{prop:recodnu2}
\partial_{\nu_{{\mathbb{B}}_n(0,r)\setminus\overline{\Omega}} 
}
Z[r_{|\partial\Omega}^t[g]]=r_{|\partial\Omega}^t[g]\qquad\forall g\in V^{-1,\alpha}(\partial\Omega)\,.
\end{equation}
Next, we choose $\varphi_r\in C^\infty({\mathbb{R}}^n)$ such that ${\mathrm{supp}}\,\varphi\subseteq {\mathbb{B}}_n(0,r)$ and $\varphi_r(x)=1$ for all $x$ in an open neighborhood $U$ of $\overline{\Omega}$ such that $\overline{U}\subseteq {\mathbb{B}}_n(0,r)$.  Next we introduce the linear operator $ \stackrel{o}{E}_{\Omega,r}$ from the subspace 
\[
X^{0,\alpha}\equiv \{f\in C^{0,\alpha}(\overline{{\mathbb{B}}_n(0,r)}\setminus\Omega)_\Delta:\,
{\mathrm{supp}}\,f\subseteq {\mathbb{B}}_n(0,r)\}
\]
of $C^{0,\alpha}(\overline{{\mathbb{B}}_n(0,r)}\setminus\Omega)_\Delta$ to $C^{0,\alpha}_{ {\mathrm{loc}} }(\overline{\Omega^-})_\Delta$ that is defined by the equality
\[
\stackrel{o}{E}_{\Omega,r}[f]\equiv\left\{
\begin{array}{ll}
 f(x) & \text{if}\ x\in \overline{{\mathbb{B}}_n(0,r)}\setminus\Omega\,,
 \\
 0 & \text{if}\ x\in {\mathbb{R}}^n\setminus\overline{{\mathbb{B}}_n(0,r)}\,,
\end{array}
\right.\qquad\forall f\in X^{0,\alpha}\,.
\]
By the definition of the topology in $C^{0,\alpha}(\overline{{\mathbb{B}}_n(0,r)}\setminus\Omega)_\Delta$ and in $C^{0,\alpha}_{ {\mathrm{loc}} }(\overline{\Omega^-})_\Delta$, the operator $ \stackrel{o}{E}_{\Omega,r}$ is continuous. Then the product Lemma \ref{lem:multc2ac0ad} implies that the map $Z_-$ from
$V^{-1,\alpha}(\partial\Omega)$ to $C^{0,\alpha}_{
{\mathrm{loc}}	}(\overline{\Omega^-})_\Delta$ that is delivered by the formula
\[
Z_-[g]\equiv\stackrel{o}{E}_{\Omega,r}\left[ \varphi_rZ[r_{|\partial\Omega}^t[g]]\right]\qquad\forall g\in V^{-1,\alpha}(\partial\Omega) \,,
\]
is linear and continuous. Then again Definition \ref{defn:endedr} implies that
\begin{equation}\label{prop:recodnu3}
 \partial_{\nu_{\Omega^-}}Z_-[g]
 =\stackrel{o}{E}_{\partial\Omega,r}^t\left[\partial_{\nu_{{\mathbb{B}}_n(0,r)\setminus\overline{\Omega}} 
}
\stackrel{o}{E}_{\Omega,r}\left[ 
\varphi_rZ[r_{|\partial\Omega}^t[g]] 
\right]\right]\qquad\forall g\in V^{-1,\alpha}(\partial\Omega)\,.
\end{equation}
Next, we prove that
\begin{equation}\label{prop:recodnu4}
\stackrel{o}{E}_{\partial\Omega,r}^t\left[\partial_{\nu_{{\mathbb{B}}_n(0,r)\setminus\overline{\Omega}} 
}
\stackrel{o}{E}_{\Omega,r}\left[ 
\varphi_rZ[r_{|\partial\Omega}^t[g]]
\right]\right]
=\stackrel{o}{E}_{\partial\Omega,r}^t\left[\partial_{\nu_{{\mathbb{B}}_n(0,r)\setminus\overline{\Omega}} 
}Z[r_{|\partial\Omega}^t[g]]\right]
\,,
\end{equation}
for all $g\in V^{-1,\alpha}(\partial\Omega)$. Let $g\in V^{-1,\alpha}(\partial\Omega)$. In order to shorten our notation, we set
$G\equiv Z[r_{|\partial\Omega}^t[g]]$. By the approximation  Lemma \ref{lem:apr0ad}, there exists a sequence $\{G_j\}_{j\in {\mathbb{N}}}$ in $C^{1,\alpha}(\overline{{\mathbb{B}}_n(0,r)}\setminus\Omega)_\Delta$ such that
 \begin{equation}\label{prop:recodnu5}
 \sup_{j\in {\mathbb{N}}}\|G_j\|_{
 C^{0,\alpha}(\overline{{\mathbb{B}}_n(0,r)}\setminus\Omega)_\Delta
 }<+\infty\,,\ \ 
 \lim_{j\to\infty}G_j=G\quad\text{in}\ C^{0,\beta}(\overline{{\mathbb{B}}_n(0,r)}\setminus\Omega)_\Delta\ \ \forall\beta\in]0,\alpha[\,.
 \end{equation}
 Then the continuity of $\partial_{\nu_{{\mathbb{B}}_n(0,r)\setminus\overline{\Omega}} 
}$ from the space $C^{0,\beta}(\overline{{\mathbb{B}}_n(0,r)}\setminus\Omega)_\Delta$ to the space $V^{-1,\beta}((\partial\Omega)\cup(\partial{\mathbb{B}}_n(0,r)))$ that is continuously embedded into the dual of the space $C^{1,\beta}((\partial\Omega)\cup(\partial{\mathbb{B}}_n(0,r)))$, the multiplication Lemma \ref{lem:multc2ac0ad} and the elementary equalities
\begin{eqnarray*}
&&\partial_{\nu_{{\mathbb{B}}_n(0,r)\setminus\overline{\Omega}} 
}
\stackrel{o}{E}_{\Omega,r}\left[ 
\varphi_rG
\right]
=\partial_{\nu_{{\mathbb{B}}_n(0,r)\setminus\overline{\Omega}} 
}
 \left[ 
\varphi_rG
\right]
\quad\text{on}\ \partial\Omega\,,
\\
&&
\partial_{\nu_{{\mathbb{B}}_n(0,r)\setminus\overline{\Omega}} 
}
 \left[\varphi_rG_j\right](x)=\partial_{\nu_{{\mathbb{B}}_n(0,r)\setminus\overline{\Omega}} 
}
 \left[G_j\right](x)\qquad\forall x\in\partial\Omega
\end{eqnarray*}
for all $j\in {\mathbb{N}}$
imply that
\begin{eqnarray*}
\lefteqn{
\langle \stackrel{o}{E}_{\partial\Omega,r}^t\left[\partial_{\nu_{{\mathbb{B}}_n(0,r)\setminus\overline{\Omega}} 
}
\stackrel{o}{E}_{\Omega,r}\left[ 
\varphi_rZ[r_{|\partial\Omega}^t[g]]
\right]\right],\psi\rangle 
=\langle   \partial_{\nu_{{\mathbb{B}}_n(0,r)\setminus\overline{\Omega}} 
}
 \left[\varphi_rG\right]
  ,\stackrel{o}{E}_{\partial\Omega,r}[\psi] \rangle 
}
\\ \nonumber
&&\qquad
=\lim_{j\to\infty}\langle \partial_{\nu_{{\mathbb{B}}_n(0,r)\setminus\overline{\Omega}} 
}
 \left[\varphi_rG_j\right] ,\stackrel{o}{E}_{\partial\Omega,r}[\psi]\rangle 
 \\ \nonumber
&&\qquad
=\lim_{j\to\infty}\int_{(\partial\Omega)\cup(\partial{\mathbb{B}}_n(0,r))}
 \partial_{\nu_{{\mathbb{B}}_n(0,r)\setminus\overline{\Omega}} 
}
 \left[\varphi_rG_j\right] \stackrel{o}{E}_{\partial\Omega,r}[\psi]\,d\sigma
 \\ \nonumber
&&\qquad
=\lim_{j\to\infty}\int_{\partial\Omega}
 \partial_{\nu_{{\mathbb{B}}_n(0,r)\setminus\overline{\Omega}} 
}
 \left[\varphi_rG_j\right]  \psi \,d\sigma
  \\ \nonumber
&&\qquad
 =\lim_{j\to\infty}\int_{\partial\Omega}
 \partial_{\nu_{{\mathbb{B}}_n(0,r)\setminus\overline{\Omega}} 
}
 \left[G_j\right]  \psi \,d\sigma
  \\ \nonumber
&&\qquad
=\lim_{j\to\infty}\int_{(\partial\Omega)\cup(\partial{\mathbb{B}}_n(0,r))}
 \partial_{\nu_{{\mathbb{B}}_n(0,r)\setminus\overline{\Omega}} 
}
 \left[G_j\right] \stackrel{o}{E}_{\partial\Omega,r}[\psi]\,d\sigma
  \\ \nonumber
&&\qquad
 =\langle  \partial_{\nu_{{\mathbb{B}}_n(0,r)\setminus\overline{\Omega}} 
}
 \left[ G\right]
  ,\stackrel{o}{E}_{\partial\Omega,r}[\psi] \rangle 
  =\langle \stackrel{o}{E}_{\partial\Omega,r}^t\left[
  \partial_{\nu_{{\mathbb{B}}_n(0,r)\setminus\overline{\Omega}} 
}
 \left[Z[r_{|\partial\Omega}^t[g]]\right]
  \right],\psi\rangle 
 \end{eqnarray*}
 for all $\psi\in C^{1,\alpha}(\partial\Omega)$ 
and thus equality (\ref{prop:recodnu4}) holds true. Then equalities (\ref{prop:recodnu2}), (\ref{prop:recodnu3}), (\ref{prop:recodnu4}) and Proposition \ref{prop:reco-1a} imply that
\begin{eqnarray*} 
\lefteqn{
 \partial_{\nu_{\Omega^-}}Z_-[g]
}
\\ \nonumber
&&\quad
=\stackrel{o}{E}_{\partial\Omega,r}^t\left[\partial_{\nu_{{\mathbb{B}}_n(0,r)\setminus\overline{\Omega}} 
}Z[r_{|\partial\Omega}^t[g]]\right]=\stackrel{o}{E}_{\partial\Omega,r}^t\left[r_{|\partial\Omega}^t[g]\right]=g
 \qquad\forall g\in V^{-1,\alpha}(\partial\Omega)
\end{eqnarray*}
 and thus equality (\ref{prop:recodnu1}) holds true. Then $Z_-$ is a right inverse of $ \partial_{\nu_{\Omega^-}}$ and accordingly $ \partial_{\nu_{\Omega^-}}$ is surjective.\hfill  $\Box$ 

\vspace{\baselineskip}

\section{A Schauder regularity theorem for   nonvariational  boundary value problems}\label{sec:nvschre}

We first introduce  a linear functional ${\mathcal{I}}_{\Omega }$ on $C^{-1,\alpha}(\overline{\Omega})$ which extends the integration in $\Omega$ to all elements of   $C^{-1,\alpha}(\overline{\Omega})$ as in \cite[Proposition~2.89]{DaLaMu21}.  
\begin{proposition}\label{prelim.wdtI}
 Let $\alpha\in]0,1]$. Let $\Omega$ be a bounded open Lipschitz subset of ${\mathbb{R}}^{n}$. Then there exists one and only one  linear and continuous  operator ${\mathcal{I}}_{\Omega}$ from the space $C^{-1,\alpha}(\overline{\Omega})$ to ${\mathbb{R}}$ such that
\begin{equation}\label{prelim.wdtI1}
{\mathcal{I}}_{\Omega}[f]=  \int_{\Omega}f_{0}\,dx+\int_{\partial\Omega}\sum_{j=1}^{n} (\nu_{\Omega})_{j}f_{j}\,d\sigma
\end{equation}
for all $f=  f_{0}+\sum_{j=1}^{n}\frac{\partial}{\partial x_{j}}f_{j}\in C^{-1,\alpha}(\overline{\Omega}) $. Moreover,
\[
{\mathcal{I}}_{\Omega}[f]=\int_{\Omega}f\,dx\qquad\forall f\in C^{0,\alpha}(\overline{\Omega})\,.
\]
\end{proposition}
Then we are ready to prove the following boundary regularity theorem.
\begin{theorem}\label{thm:schreim}
Let   $\alpha\in ]0,1[$. Let $\Omega$ be a bounded open  subset of 
 ${\mathbb{R}}^{n}$ of class $C^{1,\alpha}$. 
 Then the following statements hold. If 
 \begin{enumerate}
\item[(i)] If $u\in C^{0,\alpha}(\overline{\Omega})_\Delta$ and
\begin{equation}\label{thm:schreim1}
 \frac{\partial u}{\partial\nu_{\Omega}}\in C^{0,\alpha}(\partial\Omega)
 \end{equation}
then $u\in C^{1,\alpha}(\overline{\Omega})$.
 \item[(ii)] If $u\in C^{0,\alpha}_{{\mathrm{loc}}	}(\overline{\Omega^-})_\Delta$ and
\begin{equation}\label{thm:schreim1-}
 \frac{\partial u}{\partial\nu_{\Omega^-}}\in C^{0,\alpha}(\partial\Omega)
 \end{equation}
then $u\in C^{1,\alpha}_{{\mathrm{loc}}	}(\overline{\Omega^-})$.
\item[(iii)] If $u\in C^{0,\alpha}(\overline{\Omega})_\Delta$ and
\begin{equation}\label{thm:schreimd1}
u_{|\partial\Omega}\in C^{1,\alpha}(\partial\Omega)
 \end{equation}
then $u\in C^{1,\alpha}(\overline{\Omega})$.
 \item[(iv)] If $u\in C^{0,\alpha}_{{\mathrm{loc}}	}(\overline{\Omega^-})_\Delta$ and
\begin{equation}\label{thm:schreimd1-}
 u_{|\partial\Omega}\in C^{1,\alpha}(\partial\Omega)
 \end{equation}
then $u\in C^{1,\alpha}_{{\mathrm{loc}}	}(\overline{\Omega^-})$.
\end{enumerate}
\end{theorem}
{\bf Proof.} (i) Since $u$ solves the Neumann problem
\begin{equation}\label{thm:schreim2}
\left\{
\begin{array}{ll}
 \Delta u= \Delta u &\text{in}\ \Omega\,,
 \\
 \frac{\partial u}{\partial\nu_{\Omega}}=\frac{\partial u}{\partial\nu_{\Omega}}&\text{on}\ \partial\Omega\,,
\end{array}
\right.
 \end{equation}
 then the pair $\left(\Delta u, \frac{\partial u}{\partial\nu_{\Omega}}\right)$ satisfies the compatibility conditions
 \begin{equation}\label{thm:schreim3}
 \langle\frac{\partial u}{\partial\nu_{\Omega}},\chi_{\partial\Omega_j}\rangle=
 {\mathcal{I}}_{\Omega_j }[\Delta u]\qquad \forall j\in \{1,\dots,\kappa^+\}
  \end{equation}
  (cf. \cite[Thm.~6.4]{La24c}).  Since $\frac{\partial u}{\partial\nu_{\Omega}}\in C^{0,\alpha}(\partial\Omega)$, we have
  \[
   \langle\frac{\partial u}{\partial\nu_{\Omega}},\chi_{\partial\Omega_j}\rangle=\int_{\partial\Omega}\frac{\partial u}{\partial\nu_{\Omega}} \chi_{\partial\Omega_j}\,d\sigma\qquad\forall j\in \{1,\dots,\kappa^+\}\,.
   \]
     Then
  the solvability theorem \cite[Thm.~7.23]{DaLaMu21} for the Neumann problem in Schauder spaces implies that the Neumann problem
\begin{equation}\label{thm:schreim4}
\left\{
\begin{array}{ll}
 \Delta v= \Delta u &\text{in}\ \Omega\,,
 \\
 \frac{\partial v}{\partial\nu_{\Omega}}=\frac{\partial u}{\partial\nu_{\Omega}}&\text{on}\ \partial\Omega\,,
\end{array}
\right.
 \end{equation}
  has a solution $v\in C^{1,\alpha}(\overline{\Omega})\subseteq C^{0,\alpha}(\overline{\Omega})_\Delta$. Since both $u$ and $v$ solve the Neumann problem (\ref{thm:schreim4}) in distributional form, then \cite[Thm.~6.4]{La24c} implies that the difference  $u-v$ is constant on each connected component of $\Omega$. Hence, $u\in 
  C^{1,\alpha}(\overline{\Omega})$.

(ii) Let $r\in]0,+\infty[$ such that $\overline{\Omega}\subseteq {\mathbb{B}}_n(0,r)$. Since $u\in C^{0,\alpha}_{{\mathrm{loc}}	}(\overline{\Omega^-})_\Delta$,  Lemma \ref{lem:c0ademc1a}  implies that
$u\in C^{1,\alpha}_{{\mathrm{loc}}	}( \Omega^-)$ and accordingly
\[
{\frac{\partial u}{\partial\nu_{{\mathbb{B}}_n(0,r)\setminus\overline{\Omega}}  }}_{|
\partial {\mathbb{B}}_n(0,r)}
=
\frac{\partial u}{\partial\nu_{{\mathbb{B}}_n(0,r)}}
\in C^{0,\alpha}(\partial {\mathbb{B}}_n(0,r))\,.
\]
 Then Remark \ref{rem:localization},  Proposition \ref{prop:clnore} and assumption (\ref{thm:schreim1-}) imply that 
$\partial_{\nu_{{\mathbb{B}}_n(0,r)\setminus\overline{\Omega}}}u$ belongs to $ 
C^{0,\alpha}((\partial\Omega)\cup (\partial {\mathbb{B}}_n(0,r)))$ and thus statement (i) implies that $u_{|\overline{{\mathbb{B}}_n(0,r)}\setminus\Omega}$ belongs to $ C^{1,\alpha}(\overline{{\mathbb{B}}_n(0,r)}\setminus\Omega)$. Since the choice of $r$ is arbitrary, we conclude that  $u\in C^{1,\alpha}_{{\mathrm{loc}}	}(\overline{\Omega^-})$.

(iii) Since $u_{|\partial\Omega}\in C^{1,\alpha}(\partial\Omega)$ and $\Delta u\in C^{-1,\alpha}(\overline{\Omega})$, then the Dirichlet problem
\begin{equation}\label{thm:schreimd2}
\left\{
\begin{array}{ll}
 \Delta v= \Delta u &\text{in}\ \Omega\,,
 \\
 v=u&\text{on}\ \partial\Omega\,,
\end{array}
\right.
 \end{equation}
has a unique solution $v\in C^{1,\alpha}(\overline{\Omega})$ (cf.~\cite[Thm.~7.22]{DaLaMu21}). Since $u$ satisfies problem (\ref{thm:schreimd2}), the function $u-v$ is continuous on $\overline{\Omega}$, harmonic in $\Omega$ and vanishes on $\partial\Omega$. Then the Maximum principle implies that $u=v$. Hence, $u\in 
  C^{1,\alpha}(\overline{\Omega})$.

(iv) Since $u\in C^{0,\alpha}_{{\mathrm{loc}}	}(\overline{\Omega^-})_\Delta
$, then Lemma \ref{lem:c0ademc1a}  implies that
$u\in C^{1,\alpha}_{{\mathrm{loc}}	}( \Omega^-)$ and accordingly
$u_{|\partial{\mathbb{B}}_n(0,r)}\in C^{1,\alpha}(\partial{\mathbb{B}}_n(0,r))$.  Then   assumption (\ref{thm:schreimd1-}) implies that  
$u_{|(\partial\Omega)\cup (\partial {\mathbb{B}}_n(0,r))}\in 
C^{1,\alpha}((\partial\Omega)\cup (\partial {\mathbb{B}}_n(0,r)))$ and thus
 statement (iii) implies that $u_{|\overline{{\mathbb{B}}_n(0,r)}\setminus\Omega}\in C^{1,\alpha}(\overline{{\mathbb{B}}_n(0,r)}\setminus\Omega)$. Since the choice of $r$ is arbitrary, we conclude that $u\in C^{1,\alpha}_{{\mathrm{loc}}	}(\overline{\Omega^-})$.\hfill  $\Box$ 

\vspace{\baselineskip}

 \section{The radiation condition and uniqueness theorems for the Helmholtz equation with  Dirichlet and impedance boundary conditions}\label{sec:racouh}
Now let $k\in {\mathbb{C}}\setminus]-\infty,0]$, ${\mathrm{Im}}\, k\geq 0$. As well known in scattering theory, a function $u$ of class $C^1$ in the complement of a compact subset of ${\mathbb{R}}^{n}$ satisfies the outgoing $k$-radiation condition  provided that
\begin{equation}
\label{eq:rad1}
\lim_{x\to\infty}|x|^{\frac{n-1}{2}}(Du(x)\frac{x}{|x|}-iku(x))=0\,.
\end{equation}
It is classically known that the solutions of an exterior Dirichlet  or Neumann boundary value problem for the Helmholtz equation  that satisfy the 
 outgoing $k$-radiation condition are uniquely determined on the unbounded connected component of the exterior domain and that the same applies for the   impedance boundary conditions under appropriate assumptions on $k$ (cf., \textit{e.g.},  Colton and Kress~\cite[Thms.~3.13, 3.37]{CoKr13}, Neittaanm\"{a}ki, G.F. Roach~\cite{NiRo87}, N\'{e}d\'{e}lec~\cite[Thm.~2.6.5, p.~108]{Ne01}).  Now the classical theory requires either assumptions on the existence of the classical normal derivative on the boundary or assumptions on  the summability  of the gradient of the solution around the boundary. Here instead we are dealing with $\alpha$-H\"{o}lder continuous solutions which may not satisfy such assumptions. Thus we now prove the following counterpart of the classical uniqueness theorems. 
\begin{theorem}\label{thm:uni0ad}
 Let   $\alpha\in ]0,1[$. Let $\Omega$ be a bounded open  subset of 
 ${\mathbb{R}}^{n}$ of class $C^{1,\alpha}$. Let $k\in {\mathbb{C}}\setminus ]-\infty,0]$, ${\mathrm{Im}}\,k\geq 0$, $a\in {\mathbb{C}}$. 
 Let  $u\in C^{0,\alpha}_{	{\mathrm{loc}}	}(\overline{\Omega^-})_\Delta$ satisfy the outgoing $k$-radiation condition    and the Helmholtz equation
 \begin{equation}\label{thm:uni0ad1}
\Delta u+k^2 u=0\qquad\text{in}\ \Omega^-\,, 
\end{equation}
then the following statements hold.
\begin{enumerate}
\item[(i)] If ${\mathrm{Im}}\,(k\overline{a})\leq 0$ and $-\frac{\partial u}{\partial\nu_{\Omega^-}}+au=0$ on $ \partial\Omega$, then $u=0$ in the  unbounded  connected component of $\Omega^-$.
\item[(ii)] If $u=0$ on $ \partial\Omega$, then $u=0$ in the unbounded connected component of $\Omega^-$.
\end{enumerate}
\end{theorem}
{\bf Proof.} (i) Since $-\frac{\partial u}{\partial\nu_{\Omega^-}}=-au_{|\partial\Omega}\in C^{0,\alpha}(\partial\Omega)$, the Shauder regularization Theorem \ref{thm:schreim} (ii) implies that $u\in C^{1,\alpha}_{{\mathrm{loc}}	}(\overline{\Omega^-})$. Moreover, 
\[
{\mathrm{Im}}\,\left\{
\int_{\partial\Omega}ku\frac{\partial\overline{u}}{\partial\nu_{\Omega}}
\,d\sigma\right\}
=-{\mathrm{Im}}\,(k\overline{a})\int_{\partial\Omega}|u|^2 \,d\sigma\geq 0
\]
and accordingly the classical uniqueness Theorem \ref{thm:unext} of the Appendix implies that statement (i) holds true. 

(ii) Since $u=0\in C^{1,\alpha}(\partial\Omega)$, the Schauder regularization Theorem \ref{thm:schreim} (iv) implies that $u\in C^{1,\alpha}_{{\mathrm{loc}}	}(\overline{\Omega^-})$ and accordingly the classical uniqueness Theorem \ref{thm:unext} implies that statement (ii) holds true. \hfill  $\Box$ 

\vspace{\baselineskip}

\vspace{\baselineskip}

\section{Appendix: a classical uniqueness theorem for the Helmholtz equation}

We first introduce the following Lemma on the decay of the integral of a solution of the Helmholtz equation satisfying the outgoing radiation condition, that we prove by following   computations of  
Colton and Kress~\cite[Thm.~3.3]{CoKr13}, who considered the case in which $\Omega^-$ is connected.
\begin{proposition}
\label{prop:desq}
 Let $k\in {\mathbb{C}}\setminus \{0\}$, ${\mathrm{Im}}\,k\geq 0$. Let $\alpha\in]0,1[$. Let $\Omega$ be a bounded open Lipschitz subset of ${\mathbb{R}}^{n}$.  Let $u\in C^{1,\alpha}_{{\mathrm{loc}} }(\overline{\Omega^{-}})$ satisfy equation $\Delta u+ k^{2}u=0$ in $\Omega^{-}$ and   the outgoing $(k)$-radiation condition. Then the following statements hold
 \begin{enumerate}
 \item[(i)]
 \[
\limsup_{r\to+\infty}\int_{\partial{\mathbb{B}}_{n}(0,r)}
|u|^{2}\,d\sigma\leq
-\frac{2}{|k|^{2}}{\mathrm{Im}}\,\left(
\int_{\partial\Omega}ku\frac{\partial\overline{u}}{\partial\nu_{\Omega}}
\,d\sigma\right)\,.
 \]
\item[(ii)]
\[
\limsup_{r\to+\infty}\int_{\partial{\mathbb{B}}_{n}(0,r)}
\left|\frac{\partial u}{\partial\nu_{{\mathbb{B}}_{n}(0,r)}}\right|^{2}\,d\sigma\leq
-2{\mathrm{Im}}\,\left(
\int_{\partial\Omega}ku\frac{\partial\overline{u}}{\partial\nu_{\Omega}}
\,d\sigma\right)\,.
 \]
\item[(iii)]
\[
({\mathrm{Im}}\,k)\int_{\Omega^{-}}|u|^{2}\,d\sigma\leq
-\frac{1}{|k|^{2}}{\mathrm{Im}}\,\left(
\int_{\partial\Omega}ku\frac{\partial\overline{u}}{\partial\nu_{\Omega}}
\,d\sigma\right)\,.
\]
\item[(iv)]
\[
({\mathrm{Im}}\,k)\int_{\Omega^{-}}|Du|^{2}\,d\sigma\leq
-{\mathrm{Im}}\,\left(
\int_{\partial\Omega}ku\frac{\partial\overline{u}}{\partial\nu_{\Omega}}
\,d\sigma\right)\,.
\]
\end{enumerate}
\end{proposition}
{\bf Proof.} Let $r\in]0,+\infty[$ be such that $\overline{\Omega}\subseteq {\mathbb{B}}_{n}(0,r)$. Then we have 
\begin{eqnarray}\label{prop:desq1}
\lefteqn{
\int_{\partial{\mathbb{B}}_{n}(0,r)}\left|
\frac{\partial u}{\partial\nu_{{\mathbb{B}}_{n}(0,r)}}-iku
\right|^{2}\,d\sigma
}
\\ \nonumber
&&\ \
=
\int_{\partial{\mathbb{B}}_{n}(0,r)}\left|
\frac{\partial u}{\partial\nu_{{\mathbb{B}}_{n}(0,r)}}
\right|^{2}
+|k|^{2}|u|^{2}\,d\sigma+
2{\mathrm{Im}}\,\left(
\int_{\partial {\mathbb{B}}_{n}(0,r) }ku\frac{\partial\overline{u}}{\partial\nu_{
{\mathbb{B}}_{n}(0,r)}}
\,d\sigma\right)\,.
\end{eqnarray}
Then the first Green identity implies that
\begin{eqnarray}\label{prop:desq2}
\lefteqn{
k\int_{\partial{\mathbb{B}}_{n}(0,r)}
u\frac{\partial \overline{u}}{\partial\nu_{{\mathbb{B}}_{n}(0,r)}}\,d\sigma
-k
\int_{\partial\Omega}
u\frac{\partial \overline{u}}{\partial\nu_{\Omega}}\,d\sigma
}
\\ \nonumber
&&\qquad
=
k\int_{{\mathbb{B}}_{n}(0,r) \setminus\overline{\Omega}}|Du|^{2}\,dx
+
k\int_{{\mathbb{B}}_{n}(0,r) \setminus\overline{\Omega}}u\Delta\overline{u}\,dx
\\ \nonumber
&&\qquad
=
k\int_{{\mathbb{B}}_{n}(0,r) \setminus\overline{\Omega}}|Du|^{2}\,dx
-
\overline{k}|k|^{2}\int_{{\mathbb{B}}_{n}(0,r) \setminus\overline{\Omega}}u \overline{u}\,dx\,,
\end{eqnarray}
(cf. \textit{e.g.},~\cite[Thm.~4.2]{DaLaMu21}.) Then by taking the imaginary part in (\ref{prop:desq2}), equality
(\ref{prop:desq1}) implies that
\begin{eqnarray*}
\lefteqn{
\int_{\partial{\mathbb{B}}_{n}(0,r)}\left|
\frac{\partial u}{\partial\nu_{{\mathbb{B}}_{n}(0,r)}}-iku
\right|^{2}\,d\sigma
}
\\ \nonumber
&& 
=
\int_{\partial{\mathbb{B}}_{n}(0,r)}\left|
\frac{\partial u}{\partial\nu_{{\mathbb{B}}_{n}(0,r)}}
\right|^{2}
+|k|^{2}|u|^{2}\,d\sigma
+2{\mathrm{Im}}\,\left(
\int_{\partial\Omega}ku\frac{\partial\overline{u}}{\partial\nu_{\Omega}}
\,d\sigma\right)
\\ \nonumber
&& 
\quad+2({\mathrm{Im}}\,k)\int_{{\mathbb{B}}_{n}(0,r) \setminus{\mathrm{cl}}\Omega}|Du|^{2}\,dx
 +2({\mathrm{Im}}\,k)|k|^{2}\int_{{\mathbb{B}}_{n}(0,r) \setminus{\mathrm{cl}}\Omega}|u|^{2}\,dx\,.
\end{eqnarray*}
Hence,
\begin{eqnarray*}
\lefteqn{
\int_{\partial{\mathbb{B}}_{n}(0,r)}\left|
\frac{\partial u}{\partial\nu_{{\mathbb{B}}_{n}(0,r)}}
\right|^{2}\,d\sigma
+|k|^{2}\int_{\partial{\mathbb{B}}_{n}(0,r)}|u|^{2}\,d\sigma}
\\ \nonumber
&&
+2({\mathrm{Im}}\,k)\int_{{\mathbb{B}}_{n}(0,r) \setminus\overline{\Omega}}|Du|^{2}\,dx
 +2({\mathrm{Im}}\,k)|k|^{2}\int_{{\mathbb{B}}_{n}(0,r) \setminus\overline{\Omega}}|u|^{2}\,dx
\\ \nonumber
&&\quad
=-2{\mathrm{Im}}\,\left(
\int_{\partial\Omega}ku\frac{\partial\overline{u}}{\partial\nu_{\Omega}}
\,d\sigma\right)
+\int_{\partial{\mathbb{B}}_{n}(0,r)}\left|
\frac{\partial u}{\partial\nu_{{\mathbb{B}}_{n}(0,r)}}-iku
\right|^{2}\,d\sigma
\\ \nonumber
&&\quad
\leq-2{\mathrm{Im}}\,\left(
\int_{\partial\Omega}ku\frac{\partial\overline{u}}{\partial\nu_{\Omega}}
\,d\sigma\right)
+s_{n}
\sup_{|x|=r}
\left|
|x|^{\frac{n-1}{2}}(Du(x)\cdot\frac{x}{|x|}-ik u(x))
\right|^{2}\,.
\end{eqnarray*}
Since ${\mathrm{Im}}\,k\geq 0$, the outgoing $(k)$-radiation condition and the Monotone Convergence Theorem imply that statements (i)--(iv) hold true.\hfill  $\Box$ 

\vspace{\baselineskip}

Then one can prove the following known classical uniqueness result (see Colton and Kress~\cite[Thm.~3.12]{CoKr13} who have assumed that $\Omega^-$ is connected). For the convenience of the reader we include a  proof.
\begin{theorem}
\label{thm:unext}
Let  $k\in {\mathbb{C}}\setminus ]-\infty,0]$, ${\mathrm{Im}}\,k\geq 0$, $\alpha\in]0,1[$.  Let $\Omega$ be a bounded open subset of ${\mathbb{R}}^{n}$ of class $C^{1,\alpha}$.
If  $u\in C^{1,\alpha}_{{\mathrm{loc}} }(\overline{\Omega^{-}})$ satisfies the boundary value problem
\begin{equation}
\label{thm:unext1}
\left\{
\begin{array}{ll}
\Delta u+ k^{2}u=0 &{\text{in}}\ \Omega^{-}\,,
\\
{\mathrm{Im}}\,\left\{
\int_{\partial\Omega}ku\frac{\partial\overline{u}}{\partial\nu_{\Omega}}
\,d\sigma\right\}
 \geq 0\,,
& 
\end{array}
\right.
\end{equation}
and if $u$  satisfies the outgoing $k$-radiation condition, then $u=0$ in the unbounded connected component of $\Omega^{-}$. 
\end{theorem}
{\bf Proof.} If ${\mathrm{Im}}\,k> 0$, the statement  is an immediate consequence of Proposition \ref{prop:desq} (iii).  If ${\mathrm{Im}}\,k= 0$, Proposition \ref{prop:desq} (i) implies that
\[
\limsup_{r\to+\infty}\int_{\partial {\mathbb{B}}_{n}(0,r)}|u|^{2}\,d\sigma=0\,.
\]
Since $k\in]0,+\infty[$, Colton and Kress~\cite[Lem.~3.11]{CoKr13}
(see also Taylor~\cite[Lem.~1.2, p.~218]{Ta96}) implies that $u=0$ in a neighborhood of infinity. Since the solutions of the Helmholtz equation in $\Omega^-$ are real analytic in $\Omega^-$, then $u$ must necessarily vanish on the unbounded connected component of $\Omega^{-}$. \hfill  $\Box$ 

\vspace{\baselineskip}

  \noindent
{\bf Statements and Declarations}\\

 \noindent
{\bf Competing interests:} This paper does not have any  conflict of interest or competing interest.

 \noindent
{\bf Acknowledgement.}  

The author  acknowledges  the support of  GNAMPA-INdAM     and   of the Project funded by the European Union – Next Generation EU under the National Recovery and Resilience Plan (NRRP), Mission 4 Component 2 Investment 1.1 - Call for tender PRIN 2022 No. 104 of February, 2 2022 of Italian Ministry of University and Research; Project code: 2022SENJZ3 (subject area: PE - Physical Sciences and Engineering) ``Perturbation problems and asymptotics for elliptic differential equations: variational and potential theoretic methods''.

\end{document}